\documentclass [11pt]{amsart}

\usepackage{amsmath,amsfonts,epsfig,latexsym}
\usepackage{amssymb}
\catcode`\@=11
\newcounter{proclaim}[section]

\def\theproclaim{\thesection.\arabic{proclaim}}

\newenvironment{prop}
{\refstepcounter{proclaim}\vspace{2ex}\noindent{\bf Proposition
\theproclaim\quad}\sl}{\bigskip}

\newenvironment{thm}
{\refstepcounter{proclaim}\vspace{2ex}\noindent{\bf Theorem
\theproclaim\quad}\sl}{\bigskip}

\newenvironment{lem}
{\refstepcounter{proclaim}\vspace{2ex}\noindent{\bf Lemma
\theproclaim\quad}\sl}{\bigskip}

\newenvironment{defi}
{\refstepcounter{proclaim}\vspace{2ex}\noindent{\bf Definition
\theproclaim\quad}}{\bigskip}

\newenvironment{rem}
{\refstepcounter{proclaim}\medskip\noindent{\bf Remark
\theproclaim\quad}}{{$\Diamond$}\bigskip}

\newenvironment{pr}{\medskip\noindent{\sl Proof: }}{{\null\hfill$\Box$}\medskip}

\newenvironment{pr2}{\medskip\noindent{\sl Proof of
Theorem \ref{tono}: }}{{\null\hfill$\Box$}\medskip}

\newenvironment{pr3}{\medskip\noindent{\sl Proof of
Theorem \ref{main2}: }}{{\null\hfill$\Box$}\medskip}

\newenvironment{pr4}{\medskip\noindent{\sl Proof of
Proposition \ref{compim}: }}{{\null\hfill$\Box$}\medskip}

\newcommand{\beq}{\begin{equation}}
\newcommand{\eeq}{\end{equation}}
\newcommand{\beg}{\begin{displaymath}}
\newcommand{\eeg}{\end{displaymath}}
\newcommand{\beqa}{\begin{eqnarray*}}
\newcommand{\eeqa}{\end{eqnarray*}}

\newcommand{\W}{\mathcal{W}}

\newcommand{\1}{\mathsf{1\!\!\!\;I}}
\renewcommand{\Im}{{\rm Im}}
\renewcommand{\Re}{{\rm Re}}

\newcommand{\e}{{\rm e}}
\renewcommand{\d}{{\rm d}}
\renewcommand{\i}{{\rm i}}

\font\tenms=msbm10 scaled 1200
\font\sevenms=msbm10
\font\fivems=msbm8
\newfam\bbfam \textfont\bbfam=\tenms \scriptfont\bbfam=\sevenms
\scriptscriptfont\bbfam=\fivems

\def\1{1\!{\rm I}}
\def\t1{\tilde{\1}}

\newcommand{\R}{\mathbb{R}}
\newcommand{\C}{\mathbb{C}}
\newcommand{\N}{\mathbb{N}}
\newcommand{\Z}{\mathbb{Z}}
\newcommand{\T}{\mathbb{T}}

\def\re{\mathop{\rm Re}\nolimits}
\def\im{\mathop{\rm Im}\nolimits}

\def \coi {\mathcal{ C}^\infty}

\def\adots{\mathinner{\mkern2mu\raise 1pt\hbox{.}\mkern3mu\raise 4pt\hbox
{.}\mkern 1mu\raise 7pt\hbox{.}}}

\def\eq#1{(\ref{#1})}
\def\rf#1#2#3#4{{\parindent 3,5em
\item [{\hbox to \parindent {\enskip [#1]\hfill}}] #2: {\sl
#3}\hfill\break #4
\hfill}}

\begin{document}
\title[Semiclassical resonances for a two-level Schr\"odinger operator]
{Semiclassical resonances for a two-level Schr\"odinger operator with a conical
intersection}
\author{S.\ Fujii\'e}\thanks{Supported by MSRI and 
JSPS. KAKENHI15540149}
\address[S.\ Fujii\'e]{Mathematical Institute, Tohoku University,
Aoba, Aramaki-aza, Aoba-ku, Sendai, 980-8578 Japan}
\email{fujiie@math.tohoku.ac.jp}

\author{C.\ Lasser}\thanks{Supported by DFG priority program 1095}
\address[C.\ Lasser]{Fachbereich Mathematik, Freie Universit\"at Berlin, D-14195 Berlin, 
Germany}
\email{lasser@math.fu-berlin.de}

\author{L. N\'ed\'elec}\thanks{Supported in part by the FIM of ETHZ}
\address[L. N\'ed\'elec]{L.A.G.A., Institut Gali\-l\'ee, Univer\-sit\'e de 
Paris Nord, av. J.B. Cle\-ment,
F-93430 Villetaneuse, France IUFM de l'academie de Rouen, France}
\email{nedelec@math.univ-paris13.fr}
\maketitle

\begin{abstract}
We study the resonant set of a two-level Schr\"odinger operator 
with a linear conical intersection. This model operator can be decomposed into a direct
sum of first order systems on the real half-line. For these ordinary differential systems
we locally construct exact WKB solutions, which are connected to global solutions, amongst
which are resonant states. The main results are a generalized Bohr-Sommerfeld quantization
condition and an asymptotic description of the set of resonances as a distorted lattice.
\end{abstract}

\section{Introduction}

This paper is devoted to the semiclassical distribution of resonances of 
the two-dimen\-sional and two-level Schr\"odinger operator
$$
P = -h^2\Delta_x + V(x) = 
-h^2\Delta_x + \begin{pmatrix}x_1&x_2\\x_2&-x_1\end{pmatrix}\,.
$$

Such two-level operators appear naturally in the study of molecular spectra: 
if the positions of a molecule's nuclei are denoted by $x\in\R^n$, and $0<h\ll 1$ is 
Planck's constant devided by the square-root of the nuclear mass,  then a full
molecular Hamiltonian reads as 
$H_{\rm mol}=-h^2 \Delta_x+H_{\rm el}(x)$. For every nucleonic position $x$, 
the electronic Hamiltonian $H_{\rm el}(x)$ is an operator on the electronic degrees of
freedom.  The full operator $H_{\rm mol}$ acts on nucleonic {\em and} electronic degrees
of freedom,  that is on wave functions in $L^2(\R^N,\C)$, where $N=3\cdot(\mbox{number of
nuclei+electrons})$ is a notoriously large number.  If one considers two eigenvalues of
the electronic Hamiltonian 
$H_{\rm el}(x)$, which  are well separated from the rest of the electronic spectrum
 $\sigma(H_{\rm el}(x))$ uniformly for all $x$, then Born-Oppenheimer approximation allows
to reduce the study of the full molecular problem to the case of two-level systems acting
solely on the nucleonic degrees of freedom. The justification of this approximation can be
found in \cite{cds,kms} for the time-independent and in \cite{ha80,hajo,spte} for the
time-dependent case. For considerations from the point of view of theoretical chemistry,
we also refer to 
\cite{child}.
 
Chapter 2 in the monograph \cite{ha} gives the standard classification 
of matrix Schr\"odinger operators with eigenvalue crossings of minimal multiplicity,
 where our model operator $P$ 
appears as the normal form for a codimension two crossing: 
the real symmetric potential matrix $V(x)=V(x_1,x_2)$ depends smoothly on the two real 
parameters $x_1$ and $x_2$. It has two 
eigenvalues
$$
\pm\sqrt{x_1^2+x_2^2}=\pm|x|\,, 
$$
which coincide on the codimension two manifold $\{x=0\}\subset\R^2$. 
The graph of the mapping $x\mapsto\pm|x|$ shows two cones intersecting at the origin, 
which explains the term conical intersection, which is used in the chemical physics'
literature (see for example \cite{cfrm,dyk,ya}). Such a degeneracy, geared by two
parameters, is  generic in the sense, that it cannot be removed by symmetry preserving
perturbations. The matrix $V(x)$ is in essence Rellich's  celebrated
example of a smooth matrix, which is not smoothly diagonalizable \cite{re}.

If the crossing were not present, one might decouple the system and study 
the two scalar Hamiltonians
$$
P^\pm = -h^2\Delta_x \pm |x|\,,
$$
see \cite{m,n,b}.
$P^+$ is a Schr\"odinger operator with confining potential, and 
one has pure point spectrum only (see Appendix \ref{app:plus}). 
The lower level operator $P^-$, however, has a linearly decreasing negative 
potential, and one sees by a Mourre-type argument, that 
$P^-$ is of purely continuous spectrum. 
The full operator $P$ inherits the continuous spectrum of $P^-$, 
while spectrally echoing the discrete spectrum of $P^+$ with  resonances close to the 
real axis.

In the physical literature, zero energy wave functions close to a conical level 
crossing \cite{GA,AG} have been studied, and conical intersections have also been
addressed as generators of resonances \cite{cfrm}.

A first mathematical proof of existence of resonances for the matrix operator $P$ has
been given by one of the authors \cite{these}. She crucially used that $P$
is unitarily equivalent to the direct sum of ordinary differential 
operators
\begin{equation}
\label{dec}
\bigoplus_{\nu\in h(\Z+{1\over2})} P_{\nu},\quad
P_\nu=
\begin{pmatrix}
r^2-hD_r&\nu/r\\
\nu/r&r^2+hD_r\end{pmatrix}
\end{equation}
where we use the notation $D_r=-i\partial/\partial r$.
This decomposition is achieved by a $h$-Fourier transformation, a change 
to polar coordinates, a Pr\"ufer transformation, and a $h$-Fourier series
ansatz in the angular variable. In particular, $r$ is the length $|\xi|$ of
the dual variables $\xi=(\xi_1,\xi_2)$ of $(x_1,x_2)$ and $\nu$ is the quantum angular
momentum. 

The energy surfaces of the Hamiltonian flow associated with the eigenvalues
$r^2\pm\sqrt{\rho^2+\nu^2/r^2}$ of the symbol $p_{\nu}(r,\rho)$ of $P_{\nu}$ are
$$
\left\{(r,\rho)\in\R^+\times\R;\;\rho^2=(E-r^2)^2-\nu^2/r^2\right\}.
$$
If $\nu^2<4E^3/27$, then the function $(E-r^2)^2-\nu^2/r^2$ has
three distinct positive zeros $0<r_0<r_1<r_2$, see Appendix \ref{app:es}.
The Hamiltonian flow
consists of two components, a periodic flow passing through
$(r_0,0)$ and $(r_1,0)$ and a flow coming from infinity and
going away to infinity passing through $(r_2,0)$.
%We assume that
%$\tilde\nu=\nu/h\in\Z+1/2$ is bounded from above and
%$E$ is bounded from below by positive
%constants when $h$ tends to 0. 
Let $E>0$ be positive. 
Then $r_0$ tends to 0 while 
$r_1$ and $r_2$ tend to $\sqrt E$ as $h\to 0$: the periodic and the unbounded component
 approach each other in the semiclassical limit.
%Thus the conical
%crossing has turned into a sum of avoided crossings at
%$r=\sqrt E$.  

The resonances of the operator $P$ are defined as the eigenvalues of the complex scaled 
Hamiltonian
$$
P_\theta = 
-h^2 \e^{-2i\theta}\Delta + \e^{i\theta} V(x)\,,
$$
which is a non-selfadjoint operator with discrete spectrum independent of
the dilation parameter $\theta\in]0,{\textstyle{\pi\over3}}[$
(see for example \cite{ac,her,polo,zworski} for general theory of resonances and 
\cite{these} for the model operator $P$). In terms of the reduced operators $P_\nu$,
resonances are characterized as follows:

\begin{prop}\label{compim}
$E\in\C$ is a resonance of the
operator $P=-h^2\Delta+V$ if and only if there exist
$\nu\in h(\N - \frac 1 2)=\{\frac h2,
\frac{3h}2,\ldots\}$ and a non-trivial solution $w$ to the equation
\begin{equation}
\label{reducedequation}
P_{\nu}\, w=E\, w
\end{equation}
satisfying
\begin{equation}
\label{boundarycondition}
\lim_{r\to 0+}w(r)=0, \quad
r^2w(e^{-i\theta}r), \,\,\,w'(e^{-i\theta}r)\in L^2(\R^+,\C^2),
\end{equation}
for some $\theta\in]0,\pi/3[$.
\end{prop}

The structure of the resonant set can now be analysed by constructing a solution $u$ of the 
undilated system $P_{\nu}u= Eu$ for $E\in\C$, which vanishes at $r=0$
and is incoming at $r=+\infty$ (see Proposition \ref{quantizationcond}).

Let $S_{01}(E,h)$ be the action integral for the periodic flow
defined by
$$
S_{01}(E,h)=
\int_{r_0}^{r_1}\frac{\sqrt{\nu^2-r^2(E-r^2)^2}}rdr,
$$
where the square root is defined to be positive when $E$ is
positive. It is extended analytically into a complex
neighborhood of $\{E>0\}$, and if $|E|$ is bounded from below, the asymptotic behavior of $S_{01}(E,h)$ is given by
\begin{equation}
\label{actionasymptotic}
S_{01}(E,h)={\textstyle\frac23}i E^{3/2}+
{\textstyle\frac\pi 2}i\tilde\nu h+O(h^2|\ln h|)
\qquad(h\to0).
\end{equation}
This yields the following Bohr-Sommerfeld type
quantization condition for each angular momentum
$\tilde\nu\in\N-{1\over2}$:

\begin{thm}
\label{main1}
Let $E_0>0$ and $\tilde\nu\in\N-{1\over2}=\{\frac12, \frac32,\ldots\}$ be given.
Then there exist $\epsilon>0$, $h_0>0$ and a function $\delta(E,h): 
\{(E,h)\in\C\times\R_+;\,|E-E_0|<\epsilon,\, 0<h<h_0\}\to\C$ with $\delta(E,h)\to 0$ uniformly in $E$ as $h\to0$, such that equation
\eq{reducedequation} has a non-trivial solution satisfying
\eq{boundarycondition} if and only if $(E,h)$ satisfies
the following  quantization condition: 
\begin{equation}
\label{quantizationcondition}
\sqrt{\frac{\pi h}2}\,\tilde\nu\, e^{-i\pi/4}\, E^{-3/4}\, 
e^{2S_{01}(E,h)/h}+1=\delta(E,h).
\end{equation}
\end{thm}

Combining this result with Proposition \ref{compim} and \eq{actionasymptotic}, one obtains the following theorem about the distribution of resonances. Here, we take
$\lambda=E^{3/2}$ as spectral parameter and look for resonances in
$\{\lambda\in\C;\, a<\re \lambda<b,\,\, \im \lambda<0,\,\,\im \lambda=o(1)\,\,{\rm
as}\,\,h\to 0\}$ for arbitrarily fixed positive numbers $a,b>0$. 
For each $\tilde\nu\in\N-\frac 12$, we define
$$
\Gamma_{\tilde\nu}(h)=
\left\{\lambda\in\C;\, \lambda=\lambda_{k\tilde\nu}h
-i\,{\textstyle\frac 38}\left(h\ln{\textstyle\frac
1h}-h\ln\frac{\pi\tilde\nu^2}{2\lambda_{k\tilde\nu}h}
\right),\,
k\in\Z \,\,{\rm s.t.}\,\,a<\lambda_{k\tilde\nu}h<b\right\},
$$
where $\lambda_{k\tilde\nu}=\frac{3\pi}{16}(8k-4\tilde\nu+5)$.

\begin{thm}
\label{main2}
Let $a,b>0$.
Then, for any $N\in\N$ there exists $h_0>0$ such that for any $0<h<h_0$ and 
$\lambda\in\bigcup_{\tilde\nu\le N}\Gamma_{\tilde\nu}(h) $ there is a resonance $E$ of the operator $P$ with
$\lambda - E^{3/2} = o(h)$ uniformly for all $\lambda\in\bigcup_{\tilde\nu\le N}\Gamma_{\tilde\nu}(h) $ as $h\to0$.
\end{thm}

\begin{rem}
The integer parameters $k\in\Z$ should be large of
$O(h^{-1})$ since $a<\lambda_{k\tilde\nu}h$.
Hence, the second term of the imaginary part of $\lambda\in\Gamma_{\tilde\nu}(h)$ is of
$O(h)$ and smaller than the first term. Thus, $\Gamma_{\tilde\nu}(h)$ is an almost
horizontal sequence of complex points in the $\lambda$-plane, and 
$\bigcup_{\tilde\nu\le N}\Gamma_{\tilde\nu}(h)$ is a lattice which consists of $N$
horizontal sequences. Theorem \ref{main2} means that for a fixed positive interval $[a,b]$, we
can find as many horizontal sequences as we want for sufficiently small $h$, which are close to resonances of the operator $P$.
\end{rem}

\begin{figure}[ht]
\label{fig:res}

\vspace*{-10ex}
\includegraphics[width=0.75\textwidth]{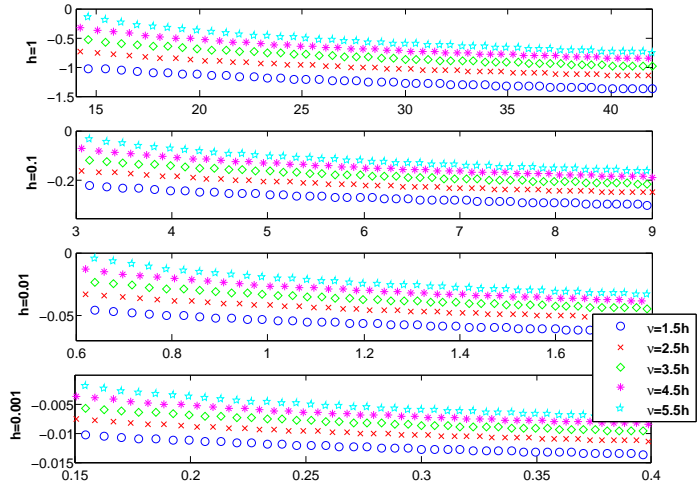}

\vspace*{-10ex}
\caption{
Resonances of the operator $P=-h^2\Delta + V(x)$. The parameter $k$ lies in
$\{11,12,\ldots,60\}$, while $\nu$ is chosen in $\{1.5h,2.5h,\ldots,5.5h\}$. The
semiclassical parameter $h$ varies from $10^{-3}$ to $1$. }
\end{figure}

The plots in Figure \ref{fig:res} illustrate the distorted lattice of
resonances given by Theorem \ref{main2}. 
The larger the angular momentum number 
$\tilde\nu\in\N-\frac{1}{2}$, the closer the resonance to the
real axis. That is, larger angular momentum numbers are
associated with longer life time of the corresponding resonant
states. Studies of the dynamical properties \cite{fermanian,ccc,cteu} of the model operator $P$ complement this observation by a quick heuristic argument:  The symbols
$p^\pm(x,\xi)=|\xi|^2\pm|x|$ of the one-level  operators
$P^\pm=-h^2\Delta\pm|x|$ induce Hamiltonian systems
$$
\dot{x}=2\xi\,,\qquad \dot{\xi}=\mp\frac{x}{|x|}  
$$
with central field. Such systems conserve angular momentum 
$x\wedge \xi=x_1\xi_2-x_2\xi_1$. Hence, only trajectories with small 
angular momentum come close to the crossing manifold $\{x=0\}$.
Moreover, angular momentum encodes closeness of trajectories to the 
crossing. 
Straightforward arguments yield that the Hamiltonian system associated with 
$p^+$ has constraint motion including periodic orbits, while the
motion corresponding to $p^-$ is unbounded.  Hence, on a heuristic
level, the classical dynamics of the decoupled systems reflect the
structure of the full operator's set of resonances: a high angular momentum
number $\widetilde\nu$ of a resonance mirrors a periodic orbit of the upper
level with high angular momentum. Such orbits in turn imply existence of
localized quasimodes and long-living resonant states. On the other hand,
small angular momentum numbers $\widetilde\nu$ correspond to orbits close
to the crossing manifold. Nearby the crossing, non-adiabatic transitions to
the unbounded motion of the minus-system are possible.
Hence, in this regime shorter life-times and resonances far away from the
real axis have to be expected. 
This heuristic point of view is in wonderful agreement with the derivation of the 
Bohr-Sommerfeld conditions (\ref{quantizationcondition}), and the prefactor 
$$
\sqrt{\frac{\pi h}2}\,\tilde\nu\, e^{-i\pi/4}\, E^{-3/4}
$$ 
can be seen as some type of Landau-Zener 
probabilty for {\em not} performing a non-adiabatic transition near the crossing. Hence, 
the conditions (\ref{quantizationcondition}) is an interpretable 
extension of quantization conditions for scalar
Schr\"{o}dinger operators to the case of a matrix-valued operator with
crossing eigenvalues. For a study relating the resonant set of the matrix operator $P$ with the spectrum of the scalar operator $P^+$ we refer to Appendix \ref{app:plus}.

The method of proving Theorem \ref{main1} and \ref{main2} is similar to \cite{rs}, while the problem resembles in some part the 
double well
problem \cite{hs}.  We construct exact solutions of the equations
$P_\nu u=Eu$, $\nu\in h(\N-\frac12)$, which are exact WKB solutions in local domains of
the complex plane. Connecting these solutions from domain to domain, we get global
solutions. More precisely, we construct an exact solution which vanishes at the origin,
and represent it, after several connection procedures, as a linear
combination of Jost solutions, which are defined at infinity. The
quantization condition of resonances will be given as the condition that
the connection coefficient
%$c^+(E,h)$ 
of the outgoing Jost solution vanishes (Proposition
\ref{quantizationcond}).

The proof of our main results, Theorem \ref{main1} and Theorem \ref{main2}, involves
several different tools and proceeds as follows:
In \S \ref{trans},  we reduce the study of the full operator $P$ to that of the ordinary 
differential operators $P_\nu$, 
$\nu\in h(\N-\frac12)$, and prove Proposition \ref{compim}.
In  \S \ref{wkb}, we extend the theory of exact WKB analysis for one-dimensional 
scalar Schr\"odinger operators used in  
\cite{GG} to a class of Schr\"odinger systems.
This exact WKB theory works in the generic situation, where the operator is without 
singularity, and turning points are "sufficiently well-separated" from each other with
respect to $h$. Our systems, however, have the origin as a regular singular point.
Moreover, the first turning point $r_0$ tends to this singularity, while the other two
turning points $r_1$ and $r_2$ approach each other at $\sqrt{E}$ as $h\to0$.
\S \ref{global} outlines the  strategy to obtain a global solution under these 
circumstances.
In \S \ref{comport}, we define Jost solutions and represent them as
exact WKB solutions.
In \S \ref{origin},
we construct an exact WKB solution vanishing at the
origin, facing the same kind of difficulty as in
\cite{langer}. The bad error estimates in Theorem
\ref{main1} and Theorem \ref{main2} come only from here.
\S \ref{norm} is devoted to the connection formula
at $\sqrt E$. 
Here, we reduce the operator 
$P_\nu$ to a well-studied microlocal normal form.
In \S \ref{distribution}, we compute the quantization condition
(Theorem \ref{main1}) from the
connection formulae obtained in the preceding sections. This condition is given in the form of Bohr-Sommerfeld using the action $S_{01}(E,h)$ between $r_0$ and $r_1$. Analysing the asymptotic behavior of the action $S_{01}(E,h)$, we finally get the semiclassical distribution of resonances 
(Theorem \ref{main2}).

\bigskip
{\sc Acknowledgments.} Part of this work has been carried out while the authors have been members of the program on semiclassical analysis at the MSRI, Berkeley. We thank Maciej Zworski for valuable  discussions.

%%%%%%%%%%%%%%%%%%%%%%%%%%%%%%%%%%%%%%%%%%%%%%%%%%%%%%%%%%%%%%%%%%%%%%%%%%%%%
\setcounter{tocdepth}{1}
{\footnotesize \tableofcontents }

\section{Reduction to the first order ordinary differential system}\indent
\label{trans}

The aim of this section is twofold: first, we define the resonances of 
\begin{equation}
\label{original}
P=-h^2\Delta+V(x),\quad 
V(x)=
\left(
\begin{array}{cc}
x_1&x_2\\
x_2&-x_1
\end{array}
\right)
\end{equation}
as the eigenvalues of its complex scaled counterpart
$$
P_\theta=-e^{-2i\,\theta}\,h^2\Delta+e^{i\,\theta}\,V(x)\,,
\qquad \theta\in ]0,{\textstyle\frac\pi3}[\,.
$$
Then, we reduce the study of the two-dimensional operators $P$ and $P_\theta$ 
to countably many one-dimensional operators $\{P_{\nu}\}_{\nu\in
h(1/2 +\Z)}$,
$$
P_{\nu}=
\left(\begin{array}{cc}
{r^2-h D_r}&{\nu/r }\\
{\nu/r }&{r^2+ h D_r}
\end{array}\right)
$$
and $\{P_{\theta,\nu}\}_{\nu\in
h(1/2 +\Z)}$
$$
P_{\theta,\nu}=
\left(\begin{array}{cc}
{e^{-2i\theta}r^2-e^{i\theta}h D_r}&{e^{i\theta}\nu/r }\\
{e^{i\theta}\nu/r }&{e^{-2i\theta}r^2+ e^{i\theta}h D_r}
\end{array}\right),
$$ 
where $D_r$ stands for $-i\partial/\partial r$. Most of the material presented here can be found in \cite{these}. We give a self-contained account of it for the convenience of the reader.

\begin{lem}
\label{lem:res}
Let $\theta\in]0,\frac\pi3[$. The operator $P_\theta$ with domain
$$
{\mathcal D}=\left\{u\in H^2(\R^2,\C^2)\,;\, |x| u\in L^2(\R^2,\C^2)\right\}
$$
is closed and has purely discrete spectrum, 
$\sigma(P_\theta)=\sigma_{\rm disc}(P_\theta)$. 
The spectrum is independent of the scaling
parameter $\theta\in ]0,\frac \pi 3[$ in the sense that 
$$
\sigma_{\rm disc}(P_\theta)=
\sigma_{\rm disc}(P_{\theta'})\qquad
(0<\theta<\theta'<{\textstyle\frac\pi3})\,.
$$
\end{lem}

\begin{pr}
We split the proof in three steps.

\medskip
{\em Closedness.}
Let $\alpha\in\C$ with $\im\alpha\neq0$ and 
$T_\alpha:=-h^2\Delta + \alpha V(x)$. Then, the following quadratic estimate holds: 
there exist positive constants $c_\alpha,b_\alpha>0$ such that
\begin{equation}
\label{eq:qe}
\|T_\alpha u\|_{L^2}^2 + c_\alpha\|u\|_{L^2}^2\ge
b_\alpha \left(\|-h^2\Delta u\|_{L^2}^2+\|\,|x|u\|_{L^2}^2\right)
\end{equation}
for all $u\in C_{\rm c}^\infty(\R^2,\C^2)$. Indeed, denoting $D=-i\nabla_x$,
$$
T_{\bar\alpha}\, T_\alpha=
(hD)^4 + \im(\alpha) \,V(hD) + 
\re(\alpha)\left(-h^2\Delta V(x)-V(x) h^2\Delta\right)+
|\alpha|^2|x|^2
\,,
$$
where one uses $i[-h^2\Delta,V(x)]=V(hD)$. Since $V(x)^2=|x|^2$, one gets
$$
0\le
\left(-h^2|\alpha|^{-\frac12}\Delta \pm |\alpha|^{\frac12}V(x)\right)^2=
|\alpha|^{-1}(hD)^4\pm (-h^2\Delta)V(x)\pm V(x)(-h^2\Delta) +|\alpha|\,|x|^2
$$
and 
$$
\re(\alpha)\left(-h^2\Delta V(x)-V(x) h^2\Delta\right)\ge-|\re(\alpha)|\,|\alpha|^{-1}\left((hD)^4+|\alpha|^2|x|^2\right)\,.
$$
Since $V(\cdot)$ behaves like $\pm|\cdot|$, one obtains 
\begin{eqnarray*}
T_{\bar\alpha}\, T_\alpha
&\ge&
(1-|\re(\alpha)|\,|\alpha|^{-1})\left((hD)^4+|\alpha|^2|x|^2\right)+
\im(\alpha) V(hD)\\
&\ge&
b_\alpha\left((hD)^4+|x|^2\right)-c_\alpha
\end{eqnarray*}
and the desired estimate (\ref{eq:qe}).
Literally the same arguments as in the proof of Theorem II.3 in \cite{her} conclude the proof, that $T_\alpha$ is a closed operator on ${\mathcal D}$.

\medskip
{\em Discrete spectrum.}
Let $\theta\in]0,\frac13[$. 
Lemma 3.1(ii) in \cite{these} proves, that the mapping $P_\theta-z_\theta:{\mathcal D}\to L^2(\R^2,\C^2)$ is bijective for all $z_\theta\in\C$ with $\im(z_\theta e^{-i\theta})>0$. By Rellich's compactness theorem, the embedding ${\mathcal D}\to L^2(\R^2,\C^2)$ is compact. Hence, 
$(P_\theta-z_\theta)^{-1}$ is compact. Let $z_\theta\in\C$ with $\im(z_\theta e^{-i\theta})>0$. For $z\in\C$ one writes
$$
(P_\theta-z)(P_\theta-z_\theta)^{-1} = 
I + (z_\theta-z)(P_\theta-z_\theta)^{-1} =: I + K_\theta(z)\,,
$$
where $\{K_\theta(z)\}_{z\in\C}$ is an analytic family of compact operators.
For $z$ sufficiently close to $z_\theta$, one has $\|K_\theta(z)\|<1$. Then, by the analytic Fredholm theorem, 
$z\mapsto(I+K_\theta(z))^{-1}$ and $z\mapsto(P_\theta-z)^{-1}$ are meromorphic in $\C$, and the residues at the poles are finite rank operators. 
Hence, $\sigma(P_\theta)=\sigma_{\rm disc}(P_\theta)$.

\medskip
{\em Independence of $\theta$.} For $\phi\in S=\left\{\phi\in\C;\; \im\phi\in]0,\frac\pi3[\right\}$ one defines
$$
H_\phi := -e^{-2\phi}h^2\Delta + e^\phi\, V(x) 
$$
as a closed operator with domain ${\mathcal D}$. 
Since $H_\phi=U_{\re\phi}P_{\theta}\,U_{-\re\phi}$ with $\theta=\im\phi$ and $U_{\re\phi}u(x)=e^{\re\phi}u(e^{\re\phi}x)$ the unitary scaling, one has $\sigma(H_\phi)=\sigma(P_\theta)$. 
Since $\{H_\phi\}_{\phi\in S}$ is a holomorphic family of type (A), 
eigenvalues $E_\phi$ of $H_\phi$ depend analytically on $\phi$, see VII.\S 3 in \cite{ka}.
By unitarity, $\phi\mapsto E_\phi$ is constant if $\im\phi$ is constant. 
Hence, the eigenvalues of $H_\phi$ and $P_\theta$ are independent of $\phi$ and $\theta$, respectively.
\end{pr}

We note, that the constants $b_\alpha,c_\alpha>0$ used in the preceding proof approach zero as $\im\alpha$ tends to zero. Hence, the above quadratic estimate does not yield closedness of the undilated operator $P=-h^2\Delta +V(x)$ on the domain ${\mathcal D}$. However, a straight forward adaption of the Faris-Lavine Theorem proves essential self-adjointness of $P$ on ${\mathcal D}$.

\begin{defi}
\label{defre}
The eigenvalues of the dilated operator $P_\theta=-e^{-2i\,\theta}\,h^2\Delta+e^{i\,\theta}\,V(x)$,
$\theta\in ]0,{\textstyle\frac\pi3}[$, with domain 
${\mathcal D}=\left\{u\in H^2(\R^2,\C^2)\,;\, |x| u\in L^2(\R^2,\C^2)\right\}$ are
called  the {\em resonances} of the operator $P=-h^2\Delta+V(x)$.
\end{defi}

The choice of dilation angle $\theta\in]0,\frac\pi 3[$ is crucial: under $x\mapsto-x$ the dilated operator $P_{\pi/3}=-e^{-2i\pi/3}h^2\Delta+e^{i\pi/3}V(x)$ is unitarily equivalent to the original operator $P$, which has no discrete spectrum, see Proposition \ref{spec}.

\begin{lem}\label{tra}
Let us fix $\theta\in ]0,\frac \pi 3[$. Then $E\in\C$ is a
resonance of $P$ if and only if there exists
$\nu\in h(\frac 1 2 +\Z)$ such that
$E\in\C$ is an eigenvalue of the operator
$P_{\theta,\nu}$
with domain
$$
{\tilde {\mathcal D}_\nu}=
\left\{w\in H^1(\R^+,\C^2);\;r^{-1}w,\, r^2 w\in L^2(\R^+,\C^2)\right\}
$$
if $\nu\neq\frac h2$, and 
$$
{\tilde {\mathcal D}_{\pm\frac h2}}=
\left\{w\in L^2(\R^+,\C^2);\;(-D_r w_1 \pm(2r)^{-1} w_2),\, 
(D_r w_2 \pm(2r)^{-1} w_1),\, r^2 w\in L^2(\R^+,\C^2)\right\}\,.
$$
\end{lem}

\begin{pr} 
Again we proceed in several steps.

\medskip{\em Fourier transformation and polar coordinates.}
Let $\hat u(\xi)\equiv ({\mathcal F}_hu)(\xi)$ be the semiclassical Fourier transform of $u$, 
$$
\hat u(\xi)=\frac 1{2\pi h}\int_{\R^2} e^{-ix\cdot\xi/h}u(x)dx.
$$
Then the dilated equation
\begin{equation}
\label{eigen1}
P_\theta u_\theta=E u_\theta
%,\quad u_\theta=U_\theta u
\end{equation}
becomes
$$
\left\{
e^{-2i\theta}|\xi|^2-
he^{i\theta}\left (
\begin{array}{cc}
D_{\xi_1} & D_{\xi_2} \\
D_{\xi_2} & -D_{\xi_1} 
\end{array}
\right )
-E
\right\}\hat u_\theta=0.
$$
Switching to polar coordinates $(\xi_1,\xi_2)=r(\cos\phi,\sin\phi)$, $r\in\R^+$, 
$\phi\in\T=\R/2\pi\Z$ one gets
$$
\left\{
e^{-2i\theta}r^2-
hA(\phi)e^{i\theta}D_r
-\frac hr A'(\phi)e^{i\theta}D_\phi-E
\right\}\hat u_\theta=0
$$
with
$$
A(\phi)=
\left (
\begin{array}{cc}
\cos\phi & \sin\phi \\
\sin\phi & -\cos\phi
\end{array}
\right ).
$$
Since the volume element $dx$ changes from $d\xi$ to $r dr d\phi$, the new Hilbert space is
$$
\left\{u\in L^2(\R^+\times\T,\C^2);\; 
\sqrt{r} u\in  L^2(\R^+\times\T,\C^2);\;
u(r,\phi+2\pi)=u(r,\phi)\right\}.
$$

\medskip{\em Half-angle rotation and conjugation by $\sqrt{r}$.}
Now put 
$$
w_\theta=\sqrt{r}\,R(-{\textstyle\frac\phi 2})\hat u_\theta,
$$
where
$$
R(\phi)=
\left (
\begin{array}{cc}
\cos\phi & -\sin\phi \\
\sin\phi & \cos\phi
\end{array}
\right )
$$
is the rotation by the angle $\phi$. 
Using the relations
$$
A(\phi )=
R(\phi )\left(\begin{array}{ll} 1 &0\\ 0& {-1}\end{array}\right) =
\left(\begin{array}{ll} 1& 0\\ 0 &{-1}\end{array}\right) R(-\phi ),
\qquad
R'(\phi)=R(\phi +{\textstyle\frac \pi 2}),
$$
we obtain the equivalent equation to (\ref{eigen1}) as
\begin{equation}
\label{final}
\left\{
e^{-2i\theta}r^2-e^{i\theta}h
\left(\begin{array}{ll} 1 &0\\ 0&
{-1}\end{array}\right)D_r-e^{i\theta}\frac hr
\left(\begin{array}{ll} 0 &1\\ 1&
{0}\end{array}\right)D_\phi-E
\right \}w_\theta=0.
\end{equation}
The new Hilbert space is
$\left\{w\in L^2(\R^+\times\T,\C^2);\; 
w(r,\phi+2\pi)=-w(r,\phi)\right\}$.

\medskip{\em Fourier series expansion.}
We expand $w_\theta$ in a Fourier series, taking into account 
that $w_\theta$ is $2\pi$ anti-periodic with respect to $\phi$:
$$
w_{\theta}(r,\phi)=\sum_{\tilde\nu\in \frac 1 2+\Z}
e^{-i\tilde \nu\phi}w_{\theta,\nu}(r).
$$
Then we see that the original eigenvalue problem (\ref{eigen1}) is reduced to that of the differential expression
$$
P_{\theta,\nu}=
e^{-2i\theta}r^2+
e^{i\theta}
\left(\begin{array}{ll} -hD_r &\frac\nu r\\ \frac\nu r&
{hD_r}\end{array}\right)\,,\qquad
\nu\in h(\Z+{\textstyle\frac12}).
$$
The resulting Hilbert space is $L^2(\R^+,\C^2)$.

\medskip
{\em Closedness of $P_{\theta,\nu}$ on $\tilde{\mathcal D}_\nu$.}
As before in the proof of Lemma \ref{lem:res} one shows a quadratic estimate, this time for the operator 
$$
T_{\alpha,\nu}:=
\begin{pmatrix}-hD_r & \frac\nu r\\ \frac\nu r & hD_r\end{pmatrix} + \alpha\, r^2
$$
with $\alpha\in\C$, $\im(\alpha)\neq0$:
there exist positive constants $b_{\alpha,\nu},c_{\alpha,\nu}>0$ such that for all $w\in C^\infty_{\rm c}(\R^+,\C^2)$
$$
\|T_{\alpha,\nu} w\|_{L^2}^2+c_{\alpha,\nu}\|w\|_{L^2}^2\ge 
\left\{
\begin{array}{ll}
b_{\alpha,\nu}
\left(
\|hD_r w\|_{L^2}^2+\|r^{-1} w\|_{L^2}^2+\| r^2 w\|_{L^2}^2 
\right)
& \mbox{if}\quad \nu\neq\pm\frac{h}{2}\,,\\*[1.5ex]
b_{\alpha,\nu}
\left(
\|\begin{pmatrix}-hD_r & \frac\nu r\\ \frac\nu r & hD_r\end{pmatrix}\! w\,\|_{L^2}^2+\| r^2 w\|_{L^2}^2 
\right)
& \mbox{if}\quad \nu=\pm\frac{h}{2}\,.
\end{array}
\right.
$$
Indeed, for $\nu\in \{\pm\frac{3h}2,\pm\frac{5h}2,\ldots\}$ one uses the lower bounds
$$
\begin{pmatrix}-hD_r & \frac\nu r\\ \frac\nu r & hD_r\end{pmatrix}^2=
-h^2\Delta_r + {\textstyle\frac{\nu^2}{r^2}} + 
ih \,{\textstyle\frac{\nu}{r^2}}
\begin{pmatrix}0 & -1\\ 1 & 0\end{pmatrix}\ge
-h^2\Delta_r + {\textstyle\frac{\nu^2-h|\nu|}{r^2}}\,,
$$
$$
\pm r^2 \begin{pmatrix}-hD_r & \frac\nu r\\ \frac\nu r & hD_r\end{pmatrix}\pm\begin{pmatrix}-hD_r & \frac\nu r\\ \frac\nu r & hD_r\end{pmatrix} r^2 \ge 
-|\alpha|^{-1}
\left(
-h^2\Delta_r + {\textstyle\frac{\nu^2-h|\nu|}{r^2}} + |\alpha|^2 r^4 
\right)\,,
$$
and
$$
i\, r^2 \begin{pmatrix}-hD_r & \frac\nu r\\ \frac\nu r & hD_r\end{pmatrix}-i \begin{pmatrix}-hD_r & \frac\nu r\\ \frac\nu r & hD_r\end{pmatrix} r^2 = 2h r
\begin{pmatrix}1 & 0\\ 0 & -1\end{pmatrix}
\ge -2rh\,.
$$
Since $\nu^2-h|\nu|>0$ for $\nu\neq\pm\frac h2$, 
one obtains for suitable $b_{\alpha,\nu},c_{\alpha,\nu}>0$
\begin{eqnarray*}
T_{\bar\alpha,\nu} T_{\alpha,\nu} &\ge& 
\left(1-|\re(\alpha)|\,|\alpha|^{-1}\right)
\left(
-h^2\Delta_r + {\textstyle\frac{\nu^2-h|\nu|}{r^2}} + |\alpha|^2 r^4 
\right)
-|\im(\alpha)| 2rh\\
&\ge&
b_{\alpha,\nu}
\left(
-h^2\Delta_r + r^{-2} + r^4 
\right)
-c_{\alpha,\nu}\,.
\end{eqnarray*}
The estimate for $\nu=\pm\frac h2$ follows analogously.
\end{pr}

\begin{pr4}
Let $E$ be a resonance of $P=-h^2\Delta+V$. By Lemma \ref{tra}, there exist $\nu\in
h(\Z+1/2)$ and $w_{\theta}\in\tilde{\mathcal D_\nu}$ such that
$P_{\theta,\nu}w_{\theta}=Ew_{\theta}$.  The origin $r=0$ is a regular singular point for
the operator $P_{\theta,\nu}$ with indicial roots $\pm\tilde\nu=\pm\nu/h$.  Hence, by the
theory of Fuchs, $w_\theta$ is a linear combination 
$$
w_\theta(r)=C_0w^0(r)+C_\infty w^\infty(r)
$$
of the two solutions $w^0(r)$ and $w^\infty(r)$ to the equation $P_{\theta,\nu}w=E w$ such that
$$
w^0(r)\sim r^{|\tilde\nu|}
\left (
\begin{array}{c}
1 \\
-{\rm sgn}(\tilde\nu)\,i
\end{array}
\right ),
\quad
w^\infty(r)\sim r^{-|\tilde\nu|}
\left (
\begin{array}{c}
1 \\
{\rm sgn}(\tilde\nu)\,i
\end{array}
\right )
\qquad
(r\to 0),
$$
see also \S 4.
Then, the condition $w_\theta\in\tilde{\mathcal D_\nu}$ implies $C_\infty=0$, and $w_\theta$ behaves
like $r^{|\tilde\nu|}$ near the origin.
In \S 5, the construction of complex WKB solutions with base 
points at infinity shows 
$$
w_\theta(r)=C_+ w^\infty_+(r)+C_- w^\infty_-(r)\,,
$$
where $w^\infty_+(r)$ is exponentially growing and $w^\infty_-(r)$ exponentially decaying as $r\to\infty$.
Hence, $C_+=0$, and $w_\theta(r)$ decays exponentially as $r\to\infty$. Setting $w(r):=w_\theta(e^{i\theta}r)$, one obtains a solution to $P_\nu w=E w$ with the claimed properties.

Let $E\in\C$ and $w$ be a solution of $P_\nu w=E w$ for some $\nu\in h(\Z+1/2)$ such that $\lim_{r\to 0+}w(r)=0$ and 
$r^2 w(e^{-i\theta}r), w'(e^{-i\theta}r)\in L^2(\R^+,\C^2)$.
The preceding arguments yield, that $w_\theta(r):=w(e^{-i\theta}r)$ is a solution to the equation $P_{\theta,\nu}w_\theta=Ew_\theta$ with 
$w_\theta\in\tilde{\mathcal D_\nu}$. By Lemma \ref{tra}, $E$ is a resonance of $P=-h^2\Delta+V$.
\end{pr4} 

\begin{rem}
\label{rem:pos}
Let $E\in\C$ and $\nu\in h(\Z+\frac{1}{2})$. Then the following equivalence holds:
$u=(u_1,u_2)$ is a solution of 
$P_\nu u=Eu$, if and only if $\widetilde u=(-u_1,u_2)$ is a solution of $P_{-\nu}\widetilde u=E\, \widetilde u$.
The same holds for the dilated operators $P_{\nu,\theta}$ and $P_{-\nu,\theta}$, $\theta\in]0,\frac{\pi}{3}[$. 
Hence, we will restrict our studies to the case $\nu\in h(\N-\frac{1}{2})$.
\end{rem}

%%%%%%%%%%%%%%%%%%%%%%%%%%%%%%%%%%%%%%%%%%%%%%%%%%%%%%%%%%%%%%%%%%%%%%%%%%%%%

\section{Exact WKB method for $2\times 2$ systems}
\label{wkb}
We now wish to find a representation formula for the solutions of $P_\nu$,
from which it is possible to deduce the asymptotic expansion in $h$. The method is
known as exact WKB method. We derive it in a somewhat more general context, and then
apply it to our specific equation.\par

We study $2\times2$ systems of first order differential equations in a complex domain $D$, which are of the form

\begin{equation}\label{eq:gen}
\left(\begin{array}{cc}
p_1(x)-{\textstyle{\frac h i}{\frac d {dx}}}&\omega(x)\\
\omega(x)&p_2(x)+{\textstyle{\frac h i}{\frac d {dx}}}
\end{array}\right)
u(x) = 0\,,
\end{equation}
or equivalently
\begin{equation}\label{eq:gen2}
{\frac h i}{\frac d {dx}}
u(x) =
\left(\begin{array}{cc}
p_1(x)&\omega(x)\\
-\omega(x)&-p_2(x)
\end{array}\right)
u(x) \,.
\end{equation}
The functions $p_1, p_2$, and $\omega$ are holomorphic in $D$.
The following considerations will lead to the construction of exact
WKB solutions for this type of systems.

%%%%%%%%%%%%%%%%%%%%%%%%%%%%%%%%

\subsection{Formal construction}
A usual change of variables (using oscillatory part) and some other basic transformations reduce the operator to a more computable one.  Let us put
$$
g_+(x)={\textstyle{\frac 1 2}}(p_1(x)+p_2(x))+\omega(x), \quad
g_-(x)=-{\textstyle{\frac 1 2}}(p_1(x)+p_2(x))+\omega(x).
$$
After conjugation by
$$
M(x)=
\exp\left({\textstyle{\frac i {2h}}}\int_{0}^x (p_1(t)-p_2(t))\,dt\right)
\left(\begin{array}{cc}1&1\\-1&1\end{array}\right)
=:m(x)
\left(\begin{array}{cc}1&1\\-1&1\end{array}\right)
$$
the system (\ref{eq:gen}) is transformed into the trace-free system
$$
{\frac h i}{\frac d {dx}}v(x) =
\left(\begin{array}{cc}
0&g_+(x)\\
-g_-(x)&0
\end{array}\right)v(x)
$$
with $u(x) = M(x)v(x)$.
Introducing a new complex coordinate
\begin{equation}\label{eq:phas}
z(x) = z(x;x_0) = \int_{x_0}^x \sqrt{
g_+(t)g_-(t)}\;dt\,,
\qquad x_0\in D\,,  \end{equation}
we look for solutions of the form
$\e^{\pm{\frac z h}}\widetilde{w}_\pm(z)$.

\begin{defi}
Let $p_1$, $p_2$, and $\omega$ be holomorphic functions in $D$.
The zeros of the function  $g_+(x)g_-(x)=-{1\over4}(p_1+p_2)^2+\omega^2$ are called the {\em
turning points} of the system (\ref{eq:gen}).
\end{defi}

We note that due to the possible presence of such turning points the  square root in the definition of $z(x)$ might be defined only {\em locally}.
By formal calculations, the amplitude vector~$\widetilde{w}_\pm(z)$ has to
satisfy
$$
{\frac h i}{d\over dz}\widetilde{w}_\pm(z) =
\left(\begin{array}{cc}
\pm i&H(z)^{-2}\\-H(z)^2&\pm i
\end{array}\right)\widetilde{w}_\pm(z)\,,
$$
where the function $H$ is given by
$$
H(z(x))=\left (\frac{g_-(x)}{g_+(x)}\right )^{1/4}\,.
$$
For a decomposition with respect to image and kernel of the preceding
system's matrix, we conjugate by
$$
P_\pm(z) \,=\, 2^{-1}
\left(\begin{array}{cc}
H(z)& \pm iH(z)^{-1}\\
H(z)&\mp  iH(z)^{-1}
\end{array}
\right)\,,\qquad
P^{-1}_\pm(z) \,=\,
\left(\begin{array}{cc}
H(z)^{-1}& H(z)^{-1}\\
\mp iH(z)&\pm  iH(z)
\end{array}
\right)
$$
and obtain a system for
$w_\pm(z) = P_\pm(z)\widetilde{w}_\pm(z)$,
$$
\frac d{dz} w_\pm(z) =
\left(\begin{array}{cc}
0&{H'(z)\over H(z)}\\
{H'(z)\over H(z)}&\mp{\textstyle{2\over h}}
\end{array}\right)w_\pm(z)\,,
$$
where $H'(z)$ is shorthand for $dH(z)/dz$.
The series ansatz
\begin{equation}\label{eq:ser}
w_\pm(z) = \sum_{n\ge0}
\left(\begin{array}{c}w_{2n,\pm}(z)\\w_{2n+1,\pm}(z)\end{array}\right)
\end{equation}
with $w_{0,\pm}\equiv1$ and for $n\ge 1$, the recurrence equations
\begin{equation}
\label{recurrence1}
\left(\frac d{dz}\pm{2\over h}\right)w_{2n+1,\pm}(z)=
{H'(z)\over H(z)}\,w_{2n,\pm}(z),
\end{equation}
\begin{equation}
\label{recurrence2}
\frac d{dz} w_{2n+2,\pm}(z)=
{H'(z)\over H(z)}\,w_{2n+1,\pm}(z)
\end{equation}
give us a formal solution up to some additive constants, which are
fixed by setting
$$
w_{n,\pm}(\widetilde{z}) = 0\,,\quad n\ge1\,,
$$
for a base point $\widetilde{z}=z(\widetilde{x})$ where $\widetilde x\in D$
is not a turning point. We note that the preceding equations for
$w_{n,\pm}$ are the same   as the ones obtained by an exact WKB
construction for scalar Schr\"odinger  equations. See for example the work
of C.~G\'erard and A.~Grigis~\cite{GG}  or T.~Ramond~\cite{R}.

Let $\Omega$ be a simply connected subset of $D$ which does not contain any
turning point. Then the function $z=z(x)$ is conformal from $\Omega$ onto
$z(\Omega)$. Assume that $\tilde z\in z(\Omega)$. If
$\Gamma_\pm(\widetilde{z},z)$ denotes a path of finite length  in
$z(\Omega)$ connecting
$\widetilde{z}$ and $z\in z(\Omega)$,
we can formally rewrite the above differential equations for $n\ge0$ as
\begin{eqnarray*}
w_{2n+1,\pm}(z)&=&
\int_{\Gamma_\pm(\widetilde{z},z)}
\exp\left(\pm{\textstyle{2\over h}}(\zeta-z)\right)
{H'(\zeta)\over H(\zeta)}\,w_{2n,\pm}(\zeta)
\,d\zeta,\\
w_{2n+2,\pm}(z)&=&
\int_{\Gamma_\pm(\widetilde{z},z)}
{H'(\zeta)\over H(\zeta)}
\,w_{2n+1,\pm}(\zeta)\,d\zeta
\end{eqnarray*}
or after iterated integration as
\begin{eqnarray*}
w_{2n+1,\pm}(z)&=&
\int_{\Gamma_\pm(\widetilde{z},z)}
\int_{\Gamma_\pm(\widetilde{z},\zeta_{2n+1})}\,\ldots\,
\int_{\Gamma_\pm(\widetilde{z},\zeta_1)}
\exp
\left(\pm{\textstyle{2\over h}(\zeta_2-\zeta_3+\ldots+\zeta_{2n+1}-z)}\right)
\times\\
&&\hspace*{24ex}
\times\;
{H'(\zeta_1)\over H(\zeta_1)}\,\ldots\,
{H'(\zeta_{2n+1})\over H(\zeta_{2n+1})}
\;d\zeta_1\ldots d\zeta_{2n+1},\\[2ex]
w_{2n+2,\pm}(z)&=&
\int_{\Gamma_\pm(\widetilde{z},z)}
\int_{\Gamma_\pm(\widetilde{z},\zeta_{2n+2})}\,\ldots\,
\int_{\Gamma_\pm(\widetilde{z},\zeta_1)}
\exp
\left(\pm{\textstyle{2\over h}}(\zeta_2-\zeta_3+\ldots-\zeta_{2n+2})\right)
\times\\
&&\hspace*{24ex}
\times\;
{H'(\zeta_1)\over H(\zeta_1)}\,\ldots\,
{H'(\zeta_{2n+2})\over H(\zeta_{2n+2})}
\;d\zeta_1\ldots d\zeta_{2n+2}\,.
\end{eqnarray*}

%%%%%%%%%%%%%%%%%%%%%%%%%%%%%%%%%%%%%%%%%%%

\subsection{Convergence, $h$-dependence, and Wronskians}
We now  give the preceding formal construction some mathematical
meaning in turning point-free compact sets $\Omega\subset D$.

\begin{lem}
For any fixed $h>0$, the formal series (\ref{eq:ser}) converges uniformly in
any compact subset of  $\Omega$, and
\begin{equation}
\label{asym}
w_\pm^{\rm even}(x,h)=\sum_{n\ge 0}w_{2n,\pm}(z(x)),\quad
w_\pm^{\rm odd}(x,h)=\sum_{n\ge 0}w_{2n+1,\pm}(z(x))
\end{equation}
are holomorphic functions in $D$.
\end{lem}

\begin{pr}
In $\Omega$, all the functions defined above are well-defined
analytic functions.
For compact subsets $K\subset\Omega$ and
$\widetilde{z}, z\in z(K)$ there exist positive constants $C_\pm^h(K)>0$
depending on the semiclassical parameter $h$ and the compact $K$ such that
\begin{equation}
\sup_{\zeta\in\Gamma_\pm(\widetilde{z},z)}
\left|\,
\exp\left(\pm{\textstyle{2\over h}}\zeta\right)
{H'(\zeta)\over H(\zeta)}
\,\right|
\;\le\; C_\pm^h(K)\,.
\label{controle}\end{equation}
If we denote the maximal length of the paths
$\Gamma_\pm(\widetilde{z},\cdot)\subset K$
in the preceding iterated integrations by $0<L<\infty$, then
$$
\sup_{z\in z(K)}|w_{n,\pm}(z)|\;\le\;
{C_\pm^h(K)^nL^n\over n!}\,,\qquad n\ge0\,,
$$
where the bound ${L^n\over n!}$ comes from the volume of a simplex with
length $L$.
\end{pr}

Thus, we have uniform convergence of the series
(\ref{eq:ser}) for $w_\pm(z)$ and exact solutions
$$
u_\pm(x) \;=\; \e^{\pm{z(x)\over h}}\,m(x)\,T_\pm(z(x))\,
\left (
\begin{array}{l}
w_\pm^{\rm even}(x) \\
w_\pm^{\rm odd}(x)
\end{array}
\right )
$$
of the original problem~(\ref{eq:gen}) on turning point free  sets~$\Omega$, where
\begin{eqnarray}
%\label{tpm}
T_\pm(z)&=&
\left(\begin{array}{cc}1&1\\-1&1\end{array}\right)
\left(\begin{array}{cc}
H(z)^{-1}& H(z)^{-1}\\
\mp iH(z)&\pm iH(z)
\end{array}
\right)\nonumber\\*[-1ex]
\label{tpm}
\\*[-1ex]
&=&
\left(\begin{array}{cc}
H(z)^{-1}\mp i H(z)&
H(z)^{-1}\pm i H(z)\\
-H(z)^{-1}\mp i H(z)&
-H(z)^{-1}\pm i H(z)
\end{array}\right)\,,\qquad z\in z(\Omega)\,.
\nonumber
\end{eqnarray}
We write these solutions $u_\pm(x)$ as 
$$
u_\pm(x;x_0,\widetilde{x})
$$ 
indicating the particular choice of the phase base point~$x_0$ in~(\ref{eq:phas}), which defines the phase function $z(x)=z(x;x_0)$, and the choice of the amplitude base point $\widetilde{z}=z(\widetilde{x})$, which is the initial point of  the path $\Gamma_\pm(\widetilde{z},\cdot)$.

For a fixed $\widetilde x\in\Omega$, let
$\Omega_\pm$ be the set of $x\in\Omega$ such that there exists a
path $\Gamma_\pm(z(\widetilde{x}),z(x))$ along which 
$x\mapsto\pm{\rm Re}\, z(x)$ increases strictly. Then,

\begin{prop}
\label{asymptotic}
The identities (\ref{asym}) for $w_\pm^{\rm even}(x,h)$ and $w_\pm^{\rm
odd}(x,h)$ give asymptotic expansions in $\Omega_\pm$. More
precisely, we have for any $\alpha\in\N$ and $N\in{\N}$
\begin{eqnarray*}
\partial^\alpha\big(
w_\pm^{\rm even}(x,h)-\sum_{n=0}^Nw_{2n,\pm}(z(x))
\big)=O(h^{N+1}),\\
\partial^\alpha\big(
w_\pm^{\rm odd}(x,h)-\sum_{n=0}^Nw_{2n+1,\pm}(z(x))
\big)=O(h^{N+2}).
\end{eqnarray*}
uniformly in compact subsets of $\Omega_\pm$. In particular,
$$
w_\pm^{\rm even}(x,h)=1+O(h),\quad
w_\pm^{\rm odd}(x,h)=O(h).
$$
\end{prop}

The proof is just the same as that of Proposition~1.2 of~\cite{GG}.
The key point is the following:
Since the iterated integrations defining $w_{n,\pm}(z)$ contain
terms of the form $\exp(\pm\zeta/h)$, one has to make sure
that $\zeta\mapsto\pm\Re(\zeta)$ is a strictly increasing function along the
path $\Gamma_\pm(\widetilde{z},z)$.
In other words, the paths $\Gamma_\pm(z(\widetilde{x}),z(x))$ have to
intersect the {\em Stokes lines}, that is the
level curves of $x\mapsto\Re(z(x))$, transversally in a suitable
direction.

\medskip
One defines the Wronskian of two $\C^2$-valued functions $u,v$ as 
${\mathcal W}(u,v)=u_1v_2-u_2v_1$.
If $z=\alpha u+\beta v$ with $\alpha,\beta\in\C$, then
$$
\alpha=\frac{\W(z,v)}{\W(u,v)}\,,\qquad
\beta=-\frac{\W(z,u)}{\W(u,v)}\,.
$$
Elementary computations give the following exact Wronskian formulas for exact WKB solutions with different phase and amplitude base points, using $w_\pm^{\rm even}$ and
$w_\pm^{\rm odd}$.

\begin{lem}
\label{lem:wron}
Let $x,x_0,y_0,\widetilde x,\widetilde y\in\Omega$. Then, 
\begin{eqnarray}
\label{wron}
\lefteqn{{\mathcal W}
(u_\pm(x;x_0,\widetilde{x}),u_\pm(x;y_0,\widetilde{y})) \;=\;
\pm 2\,i\,m(x)^2\;
\exp\left(\pm{\textstyle{1\over h}}
\left(z(x;x_0)+z(x;y_0)\right)\right)
\;}\\\nonumber
&&\hspace*{4ex}\times\left(
w^{\rm even}_{\pm}(x;x_0,\widetilde{x})
w^{\rm odd}_{\pm}(x;y_0,\widetilde{y})-
w^{\rm odd}_{\pm}(x;x_0,\widetilde{x})
w^{\rm even}_{\pm}(x;y_0,\widetilde{y})\right)\,,\\[2ex]
\label{wron1}
\lefteqn{{\mathcal W}
(u_\pm(x;x_0,\widetilde{x}),u_\mp(x;y_0,\widetilde{y})) \;=\;
\pm 2\, i\,m(x)^2\;
\exp\left(\pm{\textstyle{1\over h}}
\left(z(x;x_0)-z(x;y_0)\right)\right)
\;}\\\nonumber
&&\hspace*{4ex}\times\left(
w^{\rm even}_{\pm}(x;x_0,\widetilde{x})
w^{\rm even}_{\mp}(x;y_0,\widetilde{y})-
w^{\rm odd}_{\pm}(x;x_0,\widetilde{x})
w^{\rm odd}_{\mp}(x;y_0,\widetilde{y})\right)\,.
\end{eqnarray}
In particular, if $p_1=p_2$, then $m=1$, all the Wronskians are independent of $x$, and
we have for solutions with the same phase base point
\begin{equation}
\label{wron2}
{\mathcal W}
(u_\pm(\cdot\,;x_0,\widetilde{x}),u_\pm(\cdot\,;x_0,\widetilde{y})) \;=\;
\mp 2\,i\,
\exp\left(\pm{\textstyle{2\over h}}
z(\widetilde y;x_0)\right)
w^{\rm odd}_{\pm}(\widetilde{y}\,;x_0,\widetilde{x})\,,
\end{equation}
\begin{equation}
\label{wron3}
{\mathcal W}
(u_\pm(\cdot\,;x_0,\widetilde{x}),u_\mp(\cdot\,;x_0,\widetilde{y})) \;=\;
\pm 2\,i\;
w^{\rm even}_{\pm}(\widetilde{y}\,;x_0,\widetilde{x})\,.
\end{equation}
\end{lem}

%%%%%%%%%%%%%%%%%%%%%%%%%%%%%%%%%%%%%%%%%%%%%%%%%%%%%%%%%%%%%%%%%%%%%%%%%%%%%

\section{Formulation of the resonance condition}
\label{global}

We now go back to our system

\begin{equation}\label{eq:lmc}
(P_\nu-E)\, u(r) \;=\;
\left(\begin{array}{cc}
r^2-E-{\textstyle{h\over i}{d\over dr}}&{\nu\over r}\\
{\nu\over r}&r^2-E+{\textstyle{h\over i}{d\over dr}}
\end{array}\right)
u(r) \;=\; 0\,,\qquad r\in\R^+\,,
\end{equation}
for $E\in\C$ and $\nu\in h(\N-{1\over2})$, and study in particular
their behaviour as $r\to 0$ and $r\to\infty$.

In what follows, we use $x$ instead of $r$ and rewrite the equation
(\ref{eq:lmc}) in the form
\begin{equation}
\label{eq}
hD_xu=Au,\quad A=
\left (
\begin{array}{cc}
x^2-E & \nu/x \\
-\nu/x & E-x^2
\end{array}
\right ).
\end{equation}

The origin is a regular singular point of the equation with indicial roots
$\pm\tilde \nu$, where $\tilde \nu:=\nu/h>0$ is positive, 
see Remark \ref{rem:pos}. 
Indeed, (\ref{eq}) can be rewritten as
$$
x\frac d{dx}u=(A_0+O(x))u\quad (x\to 0),\quad A_0=
\left (
\begin{array}{cc}
0 & i\tilde\nu \\
-i\tilde\nu & 0
\end{array}
\right ),
$$
and $\pm\tilde\nu$ are the eigenvalues of $A_0$, and ${}^t(1,\mp i)$ are the
corresponding eigenvectors.
Let $u_0(x)$ be a solution corresponding to the index $\tilde\nu$,
$$
u_0(x)\sim x^{\tilde\nu}
\left (
\begin{array}{c}
1 \\
-i
\end{array}
\right )
$$
as $x\to 0$.
We will see in the next section, that there exist Jost
solutions
$f^\pm(x)$, which are characterised by their asymptotic behaviour at infinity,
\begin{equation}
\label{infinity}
f^+(x)\;\sim\;
e^{+ \,i\,(x^3-3Ex)/3h}
\left (
\begin{array}{c}
1 \\
0
\end{array}
\right ),
\quad
f^-(x)\;\sim\;
e^{- \,i\,(x^3-3Ex)/3h}
\left (
\begin{array}{c}
0 \\
1
\end{array}
\right )
\end{equation}
as $x\to+\infty$.
If $\theta\in ]0,\pi/3[$, then $f^+(xe^{-i\theta})$ is exponentially growing and 
$f^-(xe^{-i\theta})$ exponentially decaying as $x\to +\infty$. 
Since
$f^+$ and $f^-$ are linearly independent, $u_0$ can be expressed as a linear combination of these solutions,
\begin{equation}
\label{linear}
u_0(x)=c^+(E,h)f^+(x)+c^-(E,h)f^-(x).
\end{equation}
>From Definition \ref{defre} and Proposition \ref{compim}, we obtain the
following characterization of resonances:

\begin{prop}
\label{quantizationcond}
The energy $E\in\C$ is a resonance of $P$ if and only if there exists $\nu\in
h(\N-\frac12)$ with $c^+(E,h)=0$.
\end{prop}

To calculate the coefficients $c^\pm(E,h)$, which connect the
solution $u_0$ defined at the origin with the Jost solutions $f^\pm$ defined at infinity, we need some intermediate solutions, which we will construct as exact WKB solutions.
Let us recall the WKB construction of \S \ref{wkb} in this case. Note that ${\rm tr}\, A=0$ and so $m=1$. Exact WKB solutions are of the form
\begin{equation}
\label{exact}
u_\pm (x;x_0,\tilde x)=e^{\pm z(x)/h}
T_\pm(z(x))\left (
\begin{array}{l}
w_\pm^{\rm even}(x) \\
w_\pm^{\rm odd}(x)
\end{array}
\right ),
\end{equation}
where the phase function $z(x)$ is defined by
$$
z(x)=z(x;x_0)=\int_{x_0}^x\sqrt{g_+(t)g_-(t)}dt,\qquad
g_\pm(x)=\frac \nu x\mp E\pm x^2
$$
for a phase base point $x_0$, $T_\pm(z(x))$ is a $2\times 2$ matrix defined by
(\ref{tpm})
with
$$
H(z(x))=\left (\frac{g_-(x)}{g_+(x)}\right )^{1/4}=
\left (\frac{\nu +Ex- x^3}{\nu - Ex+ x^3}\right
)^{1/4},
$$
and the symbols $w_\pm^{\rm even}(x,h)=\sum_{n\ge 0}w_{2n,\pm}(z(x))$ and
$w_\pm^{\rm odd}(x,h)=\sum_{n\ge 0}w_{2n+1,\pm}(z(x))$
are constructed by recursive integrations with an amplitude  base point $\tilde x$.

There are at most six turning points, the zeros of 
$x\mapsto g_+(x)g_-(x)$, in the whole complex plane $\C_x$. They are point-symmetric with respect to the origin and denoted by
$\{\pm r_j\}_{j=0}^2$, see Appendix \ref{app:es}.
For $E>0$ fixed and $\nu>0$ sufficiently small, they are real and satisfy
$$
0<r_0<r_1<\sqrt E<r_2.
$$
As $h\to 0$, that is $\nu\to 0$, $r_0$ tends to 0 and $r_1$ and $r_2$
tend to $\sqrt E$, where $\sqrt{E}$ denotes the square root of $E\in\C$ in the right half-plane. Notice that $r_0$ and $r_1$ are zeros of $g_+(x)$ and $r_2$ is a zero of $g_-(x)$.

The Stokes curves are the level curves of the function
$x\mapsto{\rm Re}\,z(x)$.
The Stokes curves emanating from the turning points are drawn in
Figure 2.
\begin{figure}[ht]
\begin{center}
\epsfysize=5cm
\epsfbox{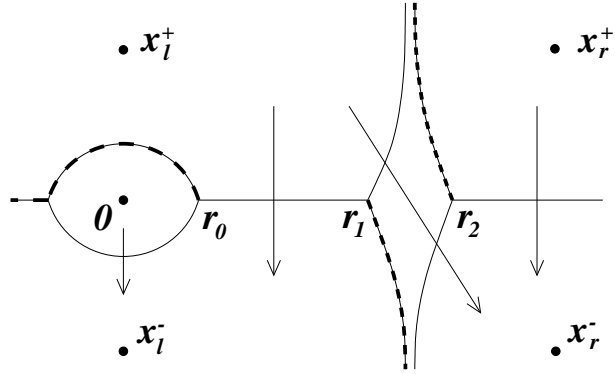}
\end{center}
\caption{The Stokes curves, that is the level curves of $x\mapsto\re z(x)$, which emanate
from the turning points $r_0$, $r_1$, $r_2$.  The arrows indicate the directions along which
$\Re\, z(x)$ increases. 
The dashed lines show the branch cuts. 
$x_l^\pm$ and $x_r^\pm$ are the amplitude base points for the constructed exact WKB
solutions.}
\end{figure}

We put branch cuts as drawn in Figure 2 for the multi-valued
functions
$$
\sqrt{g_+(x)g_-(x)}=\frac{\sqrt{\nu^2-x^2(E-x^2)^2}}x,\qquad
H(z(x))=\left (\frac{\nu+Ex-x^3}{\nu-Ex+x^3}\right )^{1/4},
$$
and suppose that
$$
\sqrt{\nu^2-x^2(E-x^2)^2}|_{x=0}=\nu,\qquad
\left (\frac{\nu+Ex-x^3}{\nu-Ex+x^3}\right )^{1/4}|_{x=0}=1.
$$
Let $E>0$ be positive, and let $\nu=\tilde \nu h>0$ be sufficiently small.
Then, for example, $\sqrt{g_+(x)g_-(x)}\in i\,\R^+$ for
$r_0<x<r_1$ and $r_2<x$, and $H(z(x))\in e^{-i\pi/4}\,\R^+$ for $r_0<x<r_1$,
while $H(z(x))\in e^{i\pi/4}\R^+$ for $r_2<x$. 

Choosing amplitude base points
$x_l^\pm$,
$x_r^\pm$ as in Figure 2 ($x_l^\pm$ are supposed to
be purely imaginary for a technical reason in the proof of
Lemma \ref{origin}), we define the following six exact WKB
solutions
$$
\begin{array}{l}
u^\pm_0(x)=u_\pm(x;r_0,x_l^\pm), \\*[0.5ex]
u_l^\pm(x)=u_\pm(x;r_1,x_l^\pm), \\*[0.5ex]
u_r^\pm(x)=u_\pm(x;r_2,x_r^\pm)
\end{array}
$$
with the turning points $r_0$, $r_1$, $r_2$ as phase base points.
The three pairs $(u_0^+(x),u_0^-(x))$, $(u_l^+(x),u_l^-(x))$,
$(u_r^+(x),u_r^-(x))$ are all linearly independent and they are
connected with $u_0$ and $(f^+,f^-)$ by transfer matrices
${}^t(c^+_0(E,h), c^-_0(E,h))$, $T_1(E,h)$,
$T_2(E,h)$, $T_3(E,h)$:
\begin{equation}
\label{con0}
u_0(x)=(u_0^+(x), u_0^-(x))
\left (
\begin{array}{c}
c^+_0(E,h) \\
c^-_0(E,h)
\end{array}
\right ),
\end{equation}
\begin{equation}
\label{con1}
(u_0^+(x), u_0^-(x))=
(u^+_l(x),u^-_l(x))T_1(E,h)
\end{equation}
\begin{equation}
\label{con2}
(u^+_l(x), u^-_l(x))=
(u^+_r(x),u^-_r(x))T_2(E,h)
\end{equation}
\begin{equation}
\label{con3}
(u^+_r(x), u^-_r(x))=
(f^+(x),f^-(x))T_3(E,h)
\end{equation}
Then, the coefficients $c^+(E,h),\,c^-(E,h)$ in \eq{linear} are given by
\begin{equation}
\label{product}
\left(\begin{array}{l}
c^+ \\
c^-
\end{array}
\right )=T_3T_2T_1
\left (
\begin{array}{l}
c^+_0 \\
c^-_0
\end{array}
\right ).
\end{equation}

Now we can describe the strategy for the remainder of this article:  The matrix $T_3$ for the transfer at infinity is computed in 
Section \ref{comport}, the connection coefficients $c^\pm_0$ at the origin in Section
\ref{origin}, and the transfer matrix
$T_2$ near $x=\sqrt E$ in Section \ref{norm}.  The matrix $T_1$ is the easiest one, and can be determined right away using Lemma \ref{lem:wron}. Because $u_0^\pm$ and $u_l^\pm$ differ only in the base point of the phase, one uses (\ref{wron}) and (\ref{wron1}) for 
$x=x_l^\pm$ to obtain
\begin{equation}
\label{T1}
T_1=
\left (
\begin{array}{cc}
e^{S_{01}/h} & 0 \\
0 & e^{-S_{01}/h}
\end{array}
\right ),
\end{equation}
where
$$
S_{01}(E,h)=\int_{r_0}^{r_1}\sqrt{g_+(x)g_-(x)}dx=
\int_{r_0}^{r_1}\frac{\sqrt{\nu^2-x^2(E-x^2)^2}}xdx
$$
is the action integral between the turning points $r_0$ and $r_1$.

%%%%%%%%%%%%%%%%%%%%%%%%%%%%%%%%%%%%%%%%%%%%%%%%%%%%%%%%%%%%%%%%%%%%%%%%%%%%

\section{Jost solutions}
\label{comport}

The Jost solutions of $P_\nu u = E u$ are characterised by their behaviour \eq{infinity} at infinity. They can be expressed as exact WKB solutions with the base points  of both
phase and amplitude placed at infinity. This fact allows us to  calculate $T_3$.

First we define the phase function with base point at infinity,
$$
z(x,\infty)=\int_{+\infty}^x
\left (\sqrt{\nu^2-t^2(E-t^2)^2}/t
-i(t^2-E)\right )dt
+{\textstyle\frac i3}(x^3-3Ex).
$$
Taking the branch of the square root into account, we see that the integral converges
absolutely, hence
$$
z(x,\infty)=\frac i3(x^3-3Ex)+o(1)\qquad (x\to +\infty).
$$

Next we define amplitudes based at infinity.
The Stokes curves are asymptotically like horizontal lines
$\{\Im x = \mbox{const.}\}$, and
$\Im x\mapsto\Re z(x)$ is a decreasing function (see Figure 2).
As in Section~3 of~\cite{R}, we choose infinite paths
$\gamma_\pm(x)$ starting from infinity and ending at $x$,
which are asymptotically like lines of the form
$\{\Im x=\mp\delta\,\Re x\}$ for some $\delta>0$,
such that $x\mapsto\mp\Re z(x)$ are strictly increasing functions
along $\gamma_\pm(x)$.
Denoting the path $z(\gamma_\pm(x))$ by $\Gamma_\pm(+\infty,z(x))$ and
setting $w_{0,\pm}\equiv1$, we inductively define $w_{n,\pm}(z)$ by
\begin{eqnarray*}
w_{2n+1,\pm}(z)&=&
\int_{\Gamma_\pm(+\infty,z)}
\exp\left(\pm{\textstyle{2\over h}}(\zeta-z)\right)
{H'(\zeta)\over H(\zeta)}\,w_{2n,\pm}(\zeta)
\,d\zeta,\\
w_{2n+2,\pm}(z)&=&
\int_{\Gamma_\pm(+\infty,z)}
{H'(\zeta)\over H(\zeta)}
\,w_{2n+1,\pm}(\zeta)\,d\zeta\,,\qquad n\ge0\,.
\end{eqnarray*}
Noticing that
$$
\frac{H'(x)}{H(x)}=
\frac{\nu}2 \,\frac{E-3x^2}{(\nu-Ex+x^3)^2}=
O(x^{-4})
$$
as $x\to+\infty$, 
one constructs well-defined exact WKB solutions $u_\infty^\pm(x)$ corresponding to these base points, proceeding as in Section 3 of \cite{R}. Up to a constant prefactor, $u_\infty^\pm(x)$ are the previously defined Jost solutions: 

\begin{lem}
\label{fandu}
Let $u_\infty^\pm$ be the exact WKB solutions with phase and amplitude base point at infinity, $f^\pm$ the Jost solutions defined in (\ref{infinity}). Then,  
$$
f^\pm(x)=\pm{\textstyle\frac 12}\,e^{\pi i/4}\,u^\pm_\infty(x)\,,
\qquad x>0\,.
$$
\end{lem}

\begin{pr}
We just check the asymptotic behavior of $u_\infty^\pm$ at infinity.
Since $H(z(x))\to e^{\pi i/4}$ as
$x\to+\infty$, we get  by an elementary calculation
$$
u^\pm_\infty(x)\;\sim\; e^{\pm i(x^3-3Ex)/3h}
\left (
\begin{array}{cc}
1 & 1 \\
-1 & 1
\end{array}
\right )
\left (
\begin{array}{cc}
e^{-\pi i/4} & e^{-\pi i/4} \\
\mp ie^{\pi i/4} & \pm ie^{\pi i/4}
\end{array}
\right )
\left (
\begin{array}{c}
1 \\
0
\end{array}
\right ),
$$
that is,
$$
u_\infty^+(x)\sim 2e^{-\pi i/4}e^{+i(x^3-3Ex)/3h}
\left (
\begin{array}{c}
1 \\ 0
\end{array}
\right ),\qquad
u_\infty^-(x)\sim -2e^{-\pi i/4}e^{-i(x^3-3Ex)/3h}
\left (
\begin{array}{c}
0 \\ 1
\end{array}
\right )
$$
as $x\to+\infty$.
\end{pr}

As an immediate consequence of the behaviour of the solutions $u^\pm_\infty(x)$ for $x\to+\infty$, we obtain that the discrete spectrum of the full operator $P$  is empty.

\begin{prop}
\label{spec}
Let $P$ be the full operator as defined in \eq{original}. Then, 
$\sigma_{\rm disc}(P)=\emptyset$. 
\end{prop}

\begin{pr}
The proof of Lemma \ref{tra} shows, that $E\in\R$ is an
eigenvalue of the self-adjoint operator $P$ if and only if there exists
$\nu\in h(\N-{1\over2})$ such that $E\in\R$ is an eigenvalue of
the ordinary differential operator $P_{\nu}$. At infinity, every distributional solution of
(\ref{eq:lmc}) has to be a linear combination of $f^\pm$.
However, for $E\in\R$ neither $f^+$ nor $f^-$ are in $L^2(\R^+,\C^2)$.
Hence, the discrete spectrum is empty.
\end{pr}

The main result of this section is the following proposition.

\begin{prop} 
\label{prop:T3}
There exists a positive $\delta>0$ independent from $E\in\C$ and $h>0$, such that the transfer matrix $T_3$ defined in (\ref{con3}) satisfies
$$
T_3(E,h)=2e^{-\pi i/4}
\left (
\begin{array}{cc}
e^{S_{2\infty}(E,h)/h}\left(1+O(h)\right)  & O(e^{-\delta/h}) \\*[1ex]
O(e^{-\delta/h}) & e^{-S_{2\infty}(E,h)/h}\left(1+O(h)\right)
\end{array}
\right )
$$
as $h\to0$, 
where $S_{2\infty}(E,h)$ is the action between $r_2$ and $+\infty$,
$$
S_{2\infty}(E,h)=
\int_{r_2}^{+\infty}\left(\sqrt{\nu^2-x^2(E-x^2)^2}/x-i(x^2-E)\right)
dx+{\textstyle\frac i3}(r_2^3-3Er_2).
$$
\end{prop}

\begin{pr}
We use the Wronskian formulas of Lemma \ref{lem:wron} and the asymptotic expansions of Proposition \ref{asymptotic} to calculate $T_3$. Setting
$$
(u_r^+, u_r^-)=:(u^+_\infty, u^-_\infty)\,\tilde T_3\,,
$$
the previous Lemma \ref{fandu} gives
$$
T_3=\left (
\begin{array}{cc}
2e^{-\pi i/4} & 0 \\
0 & -2e^{-\pi i/4}
\end{array}
\right )
\tilde T_3\,.
$$
The difference of the phases $z(x;r_2)$ and $z(x;+\infty)$ is the action $S_{2\infty}(E,h)$, and we have
$$
{\mathcal W}(u^\pm_r,u^\mp_\infty)=
\pm2i\, e^{\pm S_{2,\infty}/h}\left(1+O(h)\right),\qquad
{\mathcal W}(u^\pm_\infty,u^\mp_\infty)=
\pm2i\left(1+O(h)\right).
$$
Since there exists $\delta>0$ with
$$
{\mathcal W}(u^\pm_r,u^\pm_\infty)=O(e^{-\delta/h}),
$$
one obtains
$$
\tilde T_3=
\left (
\begin{array}{cc}
e^{S_{2\infty}/h} \left(1+O(h)\right)& O(e^{-\delta/h}) \\
O(e^{-\delta/h}) & -e^{-S_{2\infty}/h}\left(1+O(h)\right)
\end{array}
\right ).
$$
\end{pr}

%%%%%%%%%%%%%%%%%%%%%%%%%%%%%%%%%%%%%%%%%%%%%%%%%%%%%%%%%%%%%%%%%%%%%%%%%%%%%%

%%%%%%%%%%%%%%%%%%%%%%%%%%%%%%%%%%%%%%%%%%%%%%%%%%%%%%%%%%%%%%%%%%%%%%%%%%%

\section{Asymptotics of the subdominant solution near the origin}
\label{origin}

We recall that the origin is a regular singular point of $P_\nu$ with Fuchs indices $\pm\tilde\nu$, and that $u_0$ is a solution corresponding to $\tilde\nu$,
which is unique up to a constant multiple.
The purpose of this section is to calculate the connection coefficients
$c_0^\pm(E,h)$ in (\ref{con0}), i.e. to connect $u_0$ with the exact WKB
solutions $u_0^+(x)$ and $u_0^-(x)$. For this, we first need to express
$u_0$ as an exact WKB solution.

Let us look at the asymptotic behaviour of the phase 
$z(x)=z(x;r_0)$ as $x\to0$ for a fixed $h>0$. Since
$\sqrt{\nu^2-x^2(E-x^2)^2}=\nu+x\phi(x)$
with $\phi(x)$ holomorphic near $x=0$, 
$$
e^{z(x)/h}=
\exp\!\left(\int_{r_0}^x \left(\nu/t+\phi(t)\right)dt /h\right)=
C_{E,h}\,x^{\tilde\nu}\exp\!\left(\int_0^x\phi(t)dt/h\right)
$$ 
where $C_{E,h}:=r_0^{-\tilde\nu}\exp(-\int_0^{r_0}\phi(t)dt/h)>0$ is a positive constant.
Moreover $H(x)$ is a holomorphic function at $x=0$ with
$H(0)=1$. Thus, 
\begin{equation}
\label{asymp0}
e^{z(x)/h}\,T_+(z(x))\left (
\begin{array}{r}
1 \\
0
\end{array}
\right )\sim C_{E,h}\,
x^{\tilde\nu}
\left (
\begin{array}{r}
1 \\
-i
\end{array}
\right ),\qquad
x\to 0,
\end{equation}
while $h>0$ fixed.
This suggests that $u_0$ is collinear to an exact WKB solution   
of the type $+$. 
However, it is important to notice that the amplitude base point $\tilde x$ of the WKB solution  should be placed at the origin in order to have
$$
\left (
\begin{array}{r}
w^{\rm even}_+(x) \\
w^{\rm odd}_+(x)
\end{array}
\right )\sim
\left (
\begin{array}{r}
1 \\
0
\end{array}
\right ),
\qquad
x\to0,
$$
since $x\mapsto\re z(x)$ is decreasing as $x$ tends to 0 in radial directions 
(see Figure 2).
Moreover, the origin is a singular point for the equation and
we need to check that the exact WKB solution $u_+(x;r_0,0)$ is well defined, i.e.
the recurrence equations (\ref{recurrence1}), (\ref{recurrence2})
with initial value $w_{0,+}^0\equiv1$, $w_{n,+}^0(0)=0$, $n\ge 1$, define a sequence of
holomorphic functions $\{w_{n,+}^0\}_{n\ge0}$, and the series 
$\sum_{n\ge0} w_{2n,+}^0(x)$ and $\sum_{n\ge0} w_{2n+1,+}^0(x)$
converge in a neighborhood of the origin.
We rewrite (\ref{recurrence1}), (\ref{recurrence2}) with respect to $x$:
\begin{equation}
\label{recurr1}
\left(\frac d{dx}+{2\over h}\sqrt{g_+(x)g_-(x)}\right)w^0_{2n+1,+}(x)=
{H'_x(x)\over H(x)}\,w^0_{2n,+}(x),
\end{equation}
\begin{equation}
\label{recurr2}
\frac d{dx} w^0_{2n+2,+}(x)=
{H'_x(x)\over H(x)}\,w^0_{2n+1,+}(x),
\end{equation}
where $H'_x$ stands for the derivative of $H$ with respect to $x$.

These equations are of the form
$$
\frac {dw}{dx}+\frac{b(x)}xw=f(x),\quad w(0)=0
$$
with $b(x)$ and $f(x)$ given holomorphic functions at the origin.
In our case, $b\equiv0$ for (\ref{recurr2}) or
$b(x)=2x\sqrt{g_+(x)g_-(x)}/h$,
$b(0)=2\tilde\nu$ for (\ref{recurr1}).
This Cauchy problem has a unique holomorphic solution if
$\Re\, b(0)>-1$,
and the solution is given by
\begin{equation}
\label{int0}
w(x)=x\int_0^1t^{b(0)}
\exp\!\left(-x\int_t^1\tilde b(xs)ds\right)f(xt)\,dt,
\end{equation}
where $\tilde b(x)=\left(b(x)-b(0)\right)/x$.

Hence for $\tilde \nu>0$,
$\{w_{n,+}^0\}$ are uniquely determined and given by the recursive integrals
$$
w^0_{2n+2,+}=I_0(w^0_{2n+1,+}),\quad w^0_{2n+1,+}=I_1(w^0_{2n,+}),
$$
where
$$
I_0(f)=\int_0^x
\frac{H_\xi'(\xi)}{H(\xi)}f(\xi)d\xi
=\frac 12\int_0^{x}
\frac{\nu(E-3\xi^2)}{\nu^2-\xi^2(E-\xi^2)^2} f(\xi)d\xi,
$$
$$
I_1(f)=\int_0^x\exp\!\left(-\frac2h\int_\xi^x\sqrt{g_+(t)g_-(t)}dt\right)
\frac{H_\xi'(\xi)}{H(\xi)}f(\xi)d\xi\hspace{1.5cm}
$$
\begin{equation}
\label{intx}
\hspace{1.3cm}=\frac 12\int_0^{x}e^{-2\int_{\xi}^x
\sqrt{\nu^2/t^2-(E-t^2)^2}dt/h}
\frac{\nu(E-3\xi^2)}{\nu^2-\xi^2(E-\xi^2)^2} f(\xi)d\xi.
\end{equation}
It is not difficult to see that for any fixed positive $h>0$, the series
$$
w^{\rm even}_{+,0}(x,h)=\sum_{n=0}^\infty w_{2n,+}^0(x),\quad
w^{\rm odd}_{+,0}(x,h)=\sum_{n=0}^\infty w_{2n+1,+}^0(x).
$$
are absolutely convergent in a sufficiently small neighborhood
of the origin.
Hence, the function
$$
\tilde u_0(x):=e^{z(x)/h}T_+(z(x))
\left (
\begin{array}{c}
w^{\rm even}_{+,0}(x) \\*[1ex]
w^{\rm odd}_{+,0}(x)
\end{array}
\right )
$$
defines a solution to \eq{eq}.

Next, we study the asymptotic behaviour of the connection coefficients $c_0^\pm(E,h)$ as $h\to0$, using the exact WKB solution $\tilde u_0$.
For that purpose we will need some bounds on 
$w^{\rm even}_{+,0}(x,h)$ and $w^{\rm odd}_{+,0}(x,h)$. 
Since we will finally use those in the Wronskian formulas, it is enough to deal with the case when
$x$ is purely imaginary, $x=iR$ with $R>0$. We start by the following elementary estimate.

\begin{lem}
\label{element}
For any $\kappa\ge 0$, $m>0$ and $\tau>0$,  one has
$$
\int_0^\tau e^{\kappa (r-\tau)}\frac{r^{m-1}}{(1+r)^{m+1}}dr
\le \frac 1m\left (\frac\tau{1+\tau}\right )^m.
$$
\end{lem}

\begin{pr}
One integrates by parts,
\begin{eqnarray*}
\int_0^\tau e^{\kappa (r-\tau)}\frac{r^{m-1}}{(1+r)^{m+1}}dr&=&
\frac{\tau^m}{m(1+\tau)^{m+1}} -
\frac1m\int_0^\tau r^m\, 
\frac d{dr}\!\left(\frac{e^{\kappa (r-\tau)}}{(1+r)^{m+1}}\right)dr\\
&\le& \frac 1m\left (\frac\tau{1+\tau}\right )^m,
\end{eqnarray*}
using that the subtracted integral is positive.
\end{pr}

\begin{lem}
\label{estimates}
For any $E>0$, $\tau>0$, and $n\in\N$, one has
\begin{equation}
\label{ind}
|w^0_{n,+}(i\nu \tau/E)|\le \frac{K(\tau)^n}{n!}\left (
\frac{\tau}{1+\tau}\right )^{n},
\end{equation}
where
$K(\tau)=1+3\nu^2\tau^2/E^3$. 
\end{lem}

\begin{pr}
For $x=iR$, $R>0$, we have by the changes of variables $t=is$, $\xi=i\rho$,
$$
I_0(f)|_{x=iR}=\frac1{2i}\int_0^R
\frac{\nu(E+3\rho^2)}{\nu^2+\rho^2(E+\rho^2)^2}f(i\rho)d\rho\,,
$$
$$
I_1(f)|_{x=iR}=\frac1{2i}\int_0^R
\exp\left\{-\frac 2h\int_\rho^R\frac{\sqrt{\nu^2
+s^2(E+s^2)^2}}s\,ds\right\}
\frac{\nu(E+3\rho^2)}{\nu^2+\rho^2(E+\rho^2)^2}d\rho.
$$
Since
$$
\nu^2+s^2(E+s^2)^2>E^2s^2,\quad
\frac{\nu(E+3\rho^2)}{\nu^2+\rho^2(E+\rho^2)^2}\le
\frac{\nu (E+3R^2)}{\nu^2+\rho^2E^2},
$$
we have for $R=\nu\tau/E$,
\begin{equation}
\label{i0est}
|I_0(f)|_{x=i\nu\tau/E}|\le \frac{K(\tau)}2
\int_0^\tau
\frac{1}{1+r^2}\left|f\left (i\nu r/E\right )\right|dr,
\end{equation}
\begin{equation}
\label{i1est}
|I_1(f)|_{x=i\nu\tau/E}|\le \frac{K(\tau)}2
\int_0^\tau e^{2\tilde\nu (r-\tau)}
\frac{1}{1+r^2}\left|f\left (i\nu r/E\right )\right|dr.
\end{equation}
We now proceed by induction over $n\in\N$. 
Since $w_{0,+}^0\equiv1$, inequality \eq{ind} is trivially satisfied for $n=0$.
Next, we assume that \eq{ind} holds for $n=2k$.
Then, from \eq{i1est} and Lemma \ref{element}, one has
\begin{eqnarray*}
|w_{2k+1,+}^0(i\nu\tau/E)|&=&
|I_1(w_{2k,+}^0)|_{x=i\nu\tau/E}|
\le\frac{K(\tau)}2\int_0^\tau e^{2\tilde\nu(r-\tau)}\frac 1{1+r^2}
|w_{2k,+}^0(i\nu r/E)|dr\\
&\le&\frac {K(\tau)^{2k+1}}{2k!}\int_0^\tau e^{2\tilde\nu(r-\tau)}
\frac{r^{2k}}{(1+r)^{2k+2}}dr
\le \frac{K(\tau)^{2k+1}}{(2k+1)!}\left (\frac \tau{1+\tau}\right )^{2k+1}.
\end{eqnarray*}
Thus \eq{ind} holds for $n=2k+1$.
In the same way, we can show that it holds for $n=2k+2$.
\end{pr}

\begin{prop}
\label{airy}
There exists a non-zero constant $a(E,h)\neq0$ such that
$$
u_0(x)=a(E,h)\,\tilde u_0(x).
$$
The connection coefficients $c_0^\pm(E,h)$ in \eq{con0} are
analytic in $E$ near $E_0>0$ and behave as
$$
c_0^+(E,h)=a(E,h)(1+o(1)),
\quad c_0^-(E,h)=-i\,a(E,h)(1+o(1))
$$
uniformly for $E$ near $E_0>0$ as $h\to 0$. 
\end{prop}

\begin{pr}
The first part is a direct consequence of (\ref{asymp0}) and the construction of 
$\tilde u_0(x)$.

\medskip  
The second part is an Airy type connection formula at least
at the level of the principal term. Let us review briefly how to derive this
by the Wronskian formulas \eq{wron2} and \eq{wron3}.
The coefficients $c_0^+$, $c_0^-$ are
given by
$$
c_0^+= a \frac{{\mathcal W}(\tilde u_0,u_0^-)}{{\mathcal
W}(u_0^+,u_0^-)},\quad c_0^-=-a\frac{{\mathcal W}(\tilde
u_0,u_0^+)}{{\mathcal W}(u_0^+,u_0^-)}\,.
$$
Since all the solutions involved just differ in the choice of the amplitude base point, 
we can apply formulas \eq{wron2} 
and (\ref{wron3}) to get
$$
{\mathcal W}(u_0^+,u_0^-)=2i\,w_+^{\rm even}(x_l^-;r_0,x_l^+),\qquad
{\mathcal W} (\tilde u_0,u_0^-)=2i\,w_+^{\rm even}(x_l^-;r_0,0).
$$
The Wronskian ${\mathcal W}(\tilde u_0,u_0^+)$ is more delicate, since there is a branch
 cut between the origin and $x_l^+$, see Figure 2.
Hence, $u_0^+$ should be
considered on the other Riemann surface.
Let us denote by
$\hat x$ the point on the Riemann surface  continued from $x$ to
the same point turning  around $r_0$ by the angle  $-2\pi$. 
Since $g_+(x)^{1/2}=-g_+(\hat x)^{1/2}$ and $g_+(x)^{1/4}=i g_+(\hat x)^{1/4}$, we have 
$$
z(x,r_0)=-z( \hat x,r_0),\quad 
H(x)=-iH(\hat x),\quad 
T_+ (x)\,=\,i\,T_-(\hat x),
$$
and for the series summands we get  
$w_{n,\pm}(\hat x)=w_{n,\mp}(x)$.
Consequently, we have
$$
u_0^+(x^l_+,r_0,0)=iu_-(\hat x^l_+;r_0,0),
$$
which yields 
$$
{\mathcal W}(\tilde u_0,u_0^+)=-2\,w_+^{\rm even}(\hat x_l^+;r_0,0).
$$
On the other hand, we know that
$$
w_+^{\rm even}(x_l^-;r_0;x_l^+)=1+O(h),
$$
because we can take a path from $x_l^+$ to $x_l^-$, along which $x\mapsto{\rm Re}\,z(x)$ increases and which passes far
away to the right from the turning point $r_0$.
Hence
$$
c_0^+=a\,w_+^{\rm even}(x_l^-;r_0;0)\left(1+O(h)\right), \quad
c_0^-=-ia\,w_+^{\rm even}(\hat x_l^+;r_0;0)\left(1+O(h)\right),
$$ 
and it is enough for the proof to show 
\begin{equation}
\label{limit}
\lim_{h\to 0}w_+^{\rm even}(x_l^-;r_0;0)=\lim_{h\to 0}w_+^{\rm even}(\hat
x_l^+;r_0;0)=1.
\end{equation}

Let $E>0$. 
Recall that $w^{\rm even}_{+,0}=\sum_{n=0}^\infty (I_0\circ I_1)^n(1)$.
Hence if we write $x_l^+=iR$, then
\begin{eqnarray*}
w_+^{\rm even}(\hat x_l^+;r_0;0)&=&w^{\rm even}_{+,0}(iR)=
1+(I_0\circ I_1)(w^{\rm even}_{+,0})|_{x=iR}\\
&=&1+\frac
1{2iE}\int_0^\infty\chi_{(0,ER/\nu)}(r)\frac{E+3\nu^2r^2/E^2}
{1+r^2(1+\nu^2r^2/E^3)^2}
I_1(w^{\rm even}_{+,0})|_{x=i\nu r/E}\, dr,
\end{eqnarray*}
where $\chi_{(0,ER/\nu)}$ is the characteristic function of the interval
$(0,ER/\nu)$. The integrand function is dominated by an integrable function: 
indeed, Lemma \ref{estimates} gives for $r\in(0,ER/\nu)$
\begin{eqnarray*}
|w^{\rm even}_{+,0}(i\nu r/E)|&\le&\sum_{n\ge0}|w_{+,0}^{2n}(i\nu E/r)|
\le\sum_{n\ge0}\frac{K(r)^{2n}}{(2n)!}
\left (\frac{r}{1+r}\right )^{2n}\le \cosh \tilde K,\\
|I_1(w^{\rm even}_{+,0})(i\nu r/E)|&\le& \tilde K\cosh \tilde K,
\end{eqnarray*}
with $\tilde K=1+3R^2/E$, and 
$$
0\le \chi_{(0,ER/\nu)}(r)\frac{E+3\nu^2r^2/E^2}
{1+r^2(1+\nu^2r^2/E^3)^2}\le\frac{E\tilde K}{1+r^2}.
$$
On the other hand, the integrand function tends to $0$ as $h$ tends to 0
(i.e. $\nu\to 0$) for any fixed
$r>0$, since $I_1(w^{\rm even}_{+,0})|_{x=0}=0$. By Lebesgue's
dominated convergence theorem, we obtain the second identity of
\eq{limit} for $E>0$ and by analyticity for $E$ near $E_0>0$. The first identity of \eq{limit} is obtained just in the same way.
\end{pr}

%%%%%%%%%%%%%%%%%%%%%%%%%%%%%%%%%%%%%%%%%%%%%%%%%%%%%%%%%%%%%%%%%%%%%%%%%%%%%%%

\section{Connection formula near the critical point}
\label{norm}

In this section, we compute the matrix $T_2$, which connects $u^\pm_l(x)=u^\pm(x;r_1,x_l^\pm)$ with $u^\pm_r(x)=u^\pm(x;r_2,x_r^\pm)$,
$$
(u^+_l(x), u^-_l(x))=
(u^+_r(x),u^-_r(x))T_2(E,h),
$$
see \eq{con2}.
We will show

\begin{prop}
\label{7}
For all $E\in\C$ near $E_0>0$, the transfer matrix $T_2$ is of the form
\begin{equation}
\label{?}
T_2(E,h)=\left (
\begin{array}{cc}
t(E,h) & s(E,h) \\*[1ex]
-\overline{s(E,h)} & - \overline{t(E,h)}
\end{array}
\right ),
\end{equation}
and the asymptotics of $t(E,h)$ and $s(E,h)$ for $h\to0$ are given by
$$
t(E,h)=-\sqrt{\frac{\pi h}2}\,\tilde\nu\, E^{-3/4}\,e^{-i\pi/4}
+O(h|\ln h|),\qquad
s(E,h)=i+O(h).
$$
\end{prop}

The proof of Proposition \ref{7} is given in Section \ref{sec:comp}.

\begin{rem} From the symmetry properties 
\begin{equation}
\label{symmetry1}
u_r^+=i\left(\begin{array}{cc}
0&1\\1&0\end{array}\right)\overline{u_{r}^-},\qquad
u_l^-=-i\left(\begin{array}{cc}
0&1\\1&0\end{array}\right)\overline{u_{l}^+},
\end{equation}
we know {\em a priori} that $T_2$ is indeed of the claimed form (\ref{?}). Hence, it remains to prove the asymptotic behaviour of $t(E,h)$ and $s(E,h)$ as $h\to0$.
\end{rem}

The exact WKB method is not enough to compute the $h$-asymptotics of $T_2(E,h)$, because the two
turning points
$r_1$ and $r_2$ tend to the same point $x=\sqrt E$ as $h\to0$.
In other words: 
The Wronskian of two exact WKB solutions, say $u_l^+$ and $u_r^-$, is given in
terms of $w^{\rm even}_+$ computed along a path from $x_l^+$ to $x_r^-$ passing
between $r_1$ and $r_2$, but the $h$-asymptotic formula for $w^{\rm even}_+$  (Proposition \ref{asymptotic})
fails to hold because of the singularity of the function $H$ at $r_1$ and $r_2$.

Hence, one resorts to a microlocal study of the equation 
$P_\nu u=Eu$ near the point $(x,\xi)=(\sqrt E,\xi)$ for $E>0$, where 
$\xi$ is the dual variable of $x$. 
The equation $P_\nu u=Eu$ is reduced  
to a simple microlocal normal form $Qw=0$ (see Section \ref{nono}), whose solutions are well studied. From these solutions we obtain two
basis sets of microlocal solutions $(\tilde f^+, \tilde f^-)$, $(\tilde g^+, \tilde
g^-)$ of $P_\nu u=Eu$, which are related via a constant matrix $R$:
$$(\tilde f^+, \tilde f^-)R=(\tilde g^+, \tilde g^-)$$
with 
$$
R=\begin{pmatrix}p&q\\-q&-p\end{pmatrix},
$$
see Section \ref{solno}.
The exact WKB solutions $u_l^\pm$, $u_r^\pm$ are expressed in terms of
these basis sets by
$$
(u_l^+, u^-_l)=(\tilde f^+, \tilde f^-)A_l=(\tilde g^+, \tilde g^-)B_l,
$$
$$
(u_r^+, u^-_r)=(\tilde f^+, \tilde f^-)A_r=(\tilde g^+, \tilde g^-)B_r,
$$
where the constant matrices $A_{l,r}$ and $B_{l,r}$ satisfy
\begin{equation}
\label{ARB}
A_l=RB_l,\quad A_r=RB_r.
\end{equation}
Then the matrix $T_2$ is given by
\begin{equation}
\label{AB}
T_2=A_r^{-1}A_l=B_r^{-1}B_l.
\end{equation}
Hence, the $h$-asymptotics of $t(E,h)$ and $s(E,h)$ can be obtained from the study of the connection matrices $A_{l,r}$ and $B_{l,r}$  
(see Sections \ref{coni} and \ref{sec:con}).

%%%%%%%%%%%%%%%%%%%%%%%%%%%%%%%%%%%%%%%%%%%%%%%%%%%%%%%%%%%%%%%%%%%%%%%%%%%%%%%%%%%%%%%%%%%%%%%%%%%%%%%%%% 

%%%%%%%%%%%%%%%%%%%%%%%%%%%%%%%%%%%%%%%%%%%%%%%%%%%%%%%%%%%%%%%%%%%%%%%%%%%%%%%%%%%%%%%%%%%%%%%%%%%%%%%%%

\subsection{Normal form}
\label{nono}   

We  now transform the equation $(P_{\nu}-E)u=0$ near $(x,\xi)=(\sqrt{E},\xi)$, $E>0$, to a simple microlocal normal form $Qw=0$.

\begin{thm}\label{tono}
Let $E>0$ and $u(x,h)$ be a solution of $(P_\nu-E)u=0$.
Let $V$ be the metaplectic operator associated with the $\frac\pi 4$-rotation in phase space
$$
\kappa_{\frac\pi 4}:\; T^*\R\to T^*\R,\quad
(x,\xi)\mapsto{\textstyle\frac1{\sqrt2}}(x-\xi,x+\xi).
$$
There exists a locally diffeomorphic change of coordinates $x\mapsto\phi(x)=y$, $x>0$, with $\phi(\sqrt E)=0$ and  
a matrix-valued ${\mathcal C}^\infty$-symbol $M(y,h)={\rm Id}+O(h)$ such that for any cut-off function $\chi\in\coi_{\rm c}(\R)$ identically equal to $1$ in an interval around $y=0$
$$
w(y,h)= V\!\left(\chi(y) M(y,h) u(\phi^{-1}(y),h) \right)
$$
satisfies $Qw=r$, where
$$
Q=\begin{pmatrix}y & \frac\gamma{\sqrt2}\\*[1ex] -\frac{\overline\gamma}{\sqrt2} & -hD_y\end{pmatrix},
$$
$\gamma=\gamma(E,h)$ is a constant with $\gamma(E,h)=\frac{\tilde\nu}{\sqrt2} E^{-3/4} h + O(h^2)$,  
and $r(y,h)=O(h^\infty)$ uniformly in an interval around $y=0$ together with all its derivatives.
\end{thm}

\begin{rem}
The terminology {\it metaplectic operator} 
will be recalled in Appendix 
\ref{app:pr}.
\end{rem}

\begin{pr2}
We proceed in three steps to reduce the equation $P_\nu u=Eu$, that is 
\begin{equation}
\label{eq'}
hD_xu(x)=A(x)u(x),\qquad A(x)=
\left (
\begin{array}{cc}
x^2-E & \nu/x \\
-\nu/x & E-x^2
\end{array}
\right ).
\end{equation}

\medskip
{\em First step.}
One turns the quadratic diagonal entries of $A(x)$ into linear ones. 
Let $y=\phi(x)$ with
$$
\phi(x)=(x-\sqrt E) \left ({\textstyle\frac 23} (x-\sqrt E)+2\sqrt E\right
)^{1/2}.
$$
In the complex right half-plane, the function $\phi(x)$ is a biholomorphic map with 
$\phi(\sqrt E)=0$ and  
$\phi(x)\phi'(x)=x^2-E$. The function
\begin{equation}
\label{eq:psi}
\psi(y)=\psi(\phi(x))=\frac{(\frac 23(x-\sqrt E)+2\sqrt E)^{1/2}}{x(x+\sqrt E)}
\end{equation}
is analytic in a neighborhood of $y=0$, satisfying
$\psi(0)=\frac{E^{-3/4}}{\sqrt 2}$ and $\psi(\phi(x))\phi'(x)=1/x$. Moreover, 
if $u(x)$ satisfies \eq{eq'}, then $v(y)=v(\phi(x))=u(x)$ satisfies
$$
hD_yv(y) =
\left(\begin{array}{cc}
y&\nu\psi(y)\\
-\nu\psi(y)&-y
\end{array}\right)v(y)\,.
$$

\medskip
{\em Second step.}
Recall that $\nu=\tilde\nu h$ with $\tilde\nu\in\N-\frac12$.
The second step makes the off-diagonal entries constant modulo
$O(h^\infty)$ by a change of the unknown function
$$
\widetilde w(y,h)=M(y,h)v(y,h).
$$
Lemma \ref{lem:l} constructs a matrix-valued ${\mathcal C}^\infty$-symbol 
$M(y,h)={\rm Id}+O(h)$ such that $\widetilde w(y,h)$ satisfies 
\begin{equation}
\label{G}
\begin{pmatrix}hD_y-y&-\gamma\\ \overline\gamma & hD_y+y\end{pmatrix}\widetilde 
w(y,h)=r(y,h)\widetilde w(y,h)
\end{equation}
where $\gamma=\frac{\tilde\nu}{\sqrt2} E^{-3/4} h+O(h^2)$ and $r(y,h)=O(h^\infty)$ 
uniformly in an interval around $y=0$ together with all its derivatives.

\medskip
{\em Third step.}
Multiplying a cut off function $\chi$ and then operating the metaplectic operator $V$ 
from the left to equation \eq{G}, we obtain by Lemma \ref{lem:meta}
$$
Q w(y,h)=-\frac 1 {\sqrt 2}V\left(\chi(y) r(y,h) \widetilde w(y,h)-ih\chi'(y)\widetilde 
w(y,h)\right).
$$
The right hand side is of $O(h^\infty)$ uniformly in interval around 
$y=0$ together with its all derivatives. 
\end{pr2}

\begin{rem}
\label{rem:gamma}
The proof of Theorem \ref{tono} shows, that $\gamma$, $\phi$, 
and $M$ depend analytically on $E$ for $E\in\C$ near some $E_0>0$.
\end{rem}

\begin{rem}
In Appendix \ref{app:gev} we prove, that the $M_n(y)$ with 
$M(y,h)\sim\sum_{n=0}^\infty M_n(y)h^n$ 
have an analytic resummation, which is a Gevrey symbol of index $2$. Hence, the second scale $\sqrt h$, on which the two-scale Wigner measures in \cite{fermanian} are based, reappears naturally in our normal form transformation. 
\end{rem}

\subsection{Solutions of the normal form}
\label{solno}

Here we compute the solutions of the normal form.
The equation
$$
Qw\;=0,\quad w={}^t(w_1,w_2)
$$
is equivalent to
\beq\label{solnor}
w_1=-\frac {\gamma}{\sqrt{2}y}w_2\,,\quad
\frac h i y  w'_2=\frac  {{|\gamma|}^2} 2 w_2\,.\end{equation}

This is a well-studied saddle point
problem (see for example \cite{hs}, \cite{cdvp}, or section~5 in~\cite{R}).
The system \eq{solnor} has two maximal solutions
$$f^{\pm}(y) = {}^t(f_1^{\pm}(y),f_2^{\pm}(y)), $$ with
$$f_1^{\pm}(y)= -\frac {\gamma}{\sqrt{2}y}\chi_{(0,\infty)}(\pm
y)\,|y|^{\frac i {2h} |\gamma|^2 }\,,$$
$$f_2^{\pm}(y)= \chi_{(0,\infty)}(\pm y)\,|y|^{\frac i {2h} |\gamma|^2
}\,,
$$ 
where $ \chi_{(0,\infty)}$ is the characteristic function of the interval
$(0,\infty)$. Moreover, we
have two additional solutions
\begin{equation}
\label{fg}
g^{\pm}(y) =\left(\begin{array}{cc} 0&1\\1&0\end{array}\right)
{\mathcal C}{\mathcal
F}_hf^{\pm}(y),
\end{equation}
where
${\mathcal F}_hu(\eta)=\frac 1 {\sqrt{2\pi h}}\int_\R
e^{-iy\eta/h}u(y)dy$ is the $h$-Fourier
transform and ${\mathcal C}u=\bar u$ is the complex conjugate. Indeed,
we have the following identity:
$$
\left(\begin{array}{cc}0 & 1\\1&0\end{array}\right){\mathcal C}{\mathcal F}_h
Q=-Q\left(\begin{array}{cc}0 &
1\\1&0\end{array}\right){\mathcal C}{\mathcal  F}_h.
$$
Since ${\mathcal C}{\mathcal F}_h={\mathcal F}_h^{-1}{\mathcal C}$, we also
have
\begin{equation}
\label{gf}
f^{\pm}(y) =\left(\begin{array}{cc} 0&1\\1&0\end{array}\right)
{\mathcal C}{\mathcal
F}_hg^{\pm}(y).
\end{equation}
These four solutions are linearly dependent.

\begin{prop}
\label{branching}
The solutions $f^\pm$ and $g^\pm$ of $Qw=0$ are connected by
\begin{equation}
\label{scatte}
(g^+, g^-)=(f^+, f^-)R,\qquad
R=
\left(
\begin{array}{cc}
p & q \\
-q & -p
\end{array}
\right ),
\end{equation}
where
\begin{eqnarray*}
p&=&{h^{\frac 1 2-{i\over {2h}}|\gamma|^2}\over i\gamma\sqrt{\pi}}\,
\Gamma(1-{\textstyle{i\over 2h}}|\gamma|^2) \exp({\pi\over{4h}}|\gamma|^2),\\*[1ex]
q&=&{h^{\frac 1 2-{i\over {2h}}|\gamma|^2}\over i\gamma\sqrt{\pi}}\,
\Gamma(1-{\textstyle{i\over 2h}}|\gamma|^2) \exp(-{\pi\over{4h}}|\gamma|^2).
\end{eqnarray*}
\end{prop}

\begin{pr}
The proof is just the one of Proposition~5.5 in \cite{R}. 
We check $g^+_1=pf^+_1-qf^-_1$. One writes
$$
g^+_1(y)={\mathcal C}{\mathcal F_h}f_2^+(y)
=\lim_{\epsilon\to0^+} \frac1{\sqrt{2\pi h}}\int_0^\infty
e^{i (y+i\epsilon)\eta/h} \,\eta^{-\frac{i}{2h}|\gamma|^2} d\eta
$$
and substitutes $i (y+i\epsilon)\eta/h=-t$ to obtain
$$
g^+_1(y)=
\frac{h^{1/2-\frac{i}{2h}|\gamma|^2}}{\sqrt{2\pi}}
\lim_{\epsilon\to0^+} 
(\epsilon-iy)^{\frac{i}{2h}|\gamma|^2-1}
\int_{\alpha(\epsilon)}
e^{-t} \,t^{-\frac{i}{2h}|\gamma|^2} dt,
$$
where $\alpha(\epsilon)$ is the image of the interval $(0,\infty)$ under the map $\eta\mapsto (\epsilon-iy)\eta/h$. Using Euler's integral form of the gamma function $\Gamma(z)=\int_0^\infty t^{z-1} e^{-t} dt$, ${\rm Re}(z)>0$, and Cauchy's integral formula, one gets
$$
g^+_1(y)=
\frac{h^{1/2-\frac{i}{2h}|\gamma|^2}}{\sqrt{2\pi}}
(-iy)^{\frac{i}{2h}|\gamma|^2-1}\,
\Gamma(1-{\textstyle\frac{i}{2h}}|\gamma|^2)=pf^+_1(y)-qf^-_1(y).
$$
\end{pr}

\begin{rem}
In view of the relations \eq{fg} and \eq{gf}, 
the matrix $R$ satisfies
$R\bar R={\rm Id}$, that is $|p|^2-|q|^2=1$ and $p\bar q=\bar pq$.
\end{rem}

%%%%%%%%%%%%%%%%%%%%%%%%%%%%%%%%%%%%%%%%%%%%%%%%%%

%%%%%%%%%%%%%%%%%%%%%%%%%%%%%%%%%%%%%%%%%%%%%%%%%%%%%%%%%%%%%%%
\subsection{Frequency sets of the microlocal and WKB solutions}
\label{coni}

Let us study the frequency set of the microlocal solutions $f^\pm$,
$g^\pm$ and the exact WKB solutions $u_\pm^l$, $u_\pm^r$.

\medskip
First, the frequency set of the microlocal solutions are subsets of
$\R_y\times\R_\eta$ as follows:
\begin{eqnarray*}
{\rm FS}(f_1^\pm)&\subset&
\{y=0\}\cup \{\eta=0, \pm y>0\},\\
{\rm FS}(g_1^\pm)&\subset&
\{\eta=0\}\cup \{y=0, \pm\eta>0\}.
\end{eqnarray*}
Second, let $\sigma^\pm_{r,l}\subset\R_x\times\R_\xi$ be the Lagrangian manifolds defined by
\begin{equation}
\label{eq:lag}
\sigma^\pm_l=\{r_0<x<r_1,\,\,\xi=\pm(E-x^2)\},\qquad
\sigma^\pm_r=\{r_2<x,\,\,\xi=\pm(x^2-E)\}.
\end{equation}
Since $z'(x)=\sqrt{\nu^2-x^2(E-x^2)^2}/x$, one has by Lemma \ref{bkw}
$$
{\rm FS}(u^\pm_l)\cap\{r_0<x<r_1\}\subset \sigma^\pm_l \,,
\qquad
{\rm FS}(u^\pm_r)\cap\{x>r_2\}\subset \sigma^\pm_r\,.
$$

\medskip
Now we transform to the normal form of Theorem \ref{tono}, that is operate $N$ to the exact WKB solutions $u^\pm_{l,r}$,
$$
Nu(y,h)=V\!\left(\chi(y) M(y,h) u(\phi^{-1}(y),h) \right).
$$
Since $u^\pm_{l,r}$ are solutions of the equation $(P_\nu
-E)u=0$, we see that for any $k\ge 0$
$$
D_y^kQ(Nu^\pm_{l,r})=O(h^\infty).
$$
Hence, $Nu^\pm_{l,r}$ are microlocal solutions of $Qw=0$ near $(0,0)$.
Since the vector space of such solutions is two-dimensional, see Proposition 17 of \cite{cdvp}, there exist 
matrix-valued $\coi$-symbols 
\begin{equation}
\label{eq:mat}
A_l=(\alpha_{jk}^l),\quad 
A_r=(\alpha_{jk}^r),\quad
B_l=(\beta_{jk}^l),\quad 
B_r=(\beta_{jk}^r),
\end{equation}
such that microlocally near $(0,0)$
\begin{eqnarray*}
(Nu^+_l, Nu^-_l)
&=&(f^+, f^-)A_l\;=\;(g^+, g^-)B_l,\\
(Nu^+_r, Nu^-_r)
&=&(f^+, f^-)A_r\;=\;(g^+, g^-)B_r.
\end{eqnarray*}
Returning back to the $(x,\xi)$ variables, i.e. operating $N^{-1}$ from
the left, we have microlocally near
$(\sqrt E,0)$,
\begin{eqnarray*}
(u^+_l, u^-_l)
&=&(\tilde f^+, \tilde f^-)A_l\;=\;(\tilde g^+, \tilde g^-)B_l,\\
(u^+_r, u^-_r)
&=&(\tilde f^+, \tilde f^-)A_r\;=\;(\tilde g^+, \tilde g^-)B_r,
\end{eqnarray*}
where
\begin{equation}
\label{eq:tilde}
\tilde f^\pm=N^{-1}(\tilde\chi f^\pm),\quad 
\tilde g^\pm=N^{-1}(\tilde\chi g^\pm)
\end{equation}
and $\tilde\chi\in\coi_{\rm c}(\R)$ a cut-off function, which is  identically equal to $1$ near $y=0$ and satisfies ${\rm supp}(\tilde\chi)\subset{\rm supp}(\chi)$.

\medskip
By Lemma \ref{lem:meta}, we have 
$$
{\rm FS}(\tilde f^\pm)=
\kappa_\phi^{-1}\big(\kappa_{\frac\pi4}{\rm FS}(\tilde\chi f^\pm)\big),\quad
{\rm FS}(\tilde g^\pm)=
\kappa_\phi^{-1}\big(\kappa_{\frac\pi4}{\rm FS}(\tilde\chi g^\pm)\big),
$$
where 
$$
\kappa_\phi^{-1}:\; T^*\R\to T^*\R,\quad (x,\xi)\mapsto
\left(\phi^{-1}(x),\xi\,\phi'(\phi^{-1}(x))\right)
$$
is the inverse of the canonical transformation $\kappa_\phi(x,\xi)=(\phi(x),\xi/\phi'(x))$ 
associated with $\phi$. Clearly,
\begin{eqnarray*}
\kappa_{\frac\pi4}\,{\rm FS}(f_1^\pm)&\subset&
\{\eta=-y\}\cup\{\eta=y,\, \pm y>0\},\\
\kappa_{\frac\pi4}\,{\rm FS}(g_1^\pm)&\subset&
\{\eta=y\}\cup \{\eta=-y,\,\mp y>0\}.
\end{eqnarray*}
Since $\phi(x)\phi'(x)=x^2-E$, this yields
\begin{eqnarray*}
{\rm FS}(\tilde f^+)\cap U
\subset\sigma^-_l\cup\sigma^+_r\cup\sigma^-_r,&&\qquad
{\rm FS}(\tilde f^-)\cap U
\subset\sigma^+_l\cup\sigma^-_l\cup\sigma^+_r,\\
{\rm FS}(\tilde g^+)\cap U
\subset\sigma^+_l\cup\sigma^+_r\cup\sigma^-_r,&&\qquad
{\rm FS}(\tilde g^-)\cap U
\subset\sigma^+_l\cup\sigma^-_l\cup\sigma^-_r
\end{eqnarray*}
with $U=\{r_0<x<r_1\;\mbox{or}\; x>r_2\}$, see Figure 3.
$\sigma_l^+$ is neither contained in ${\rm FS}(u_l^-)$ nor in 
${\rm FS}(\tilde f^+)$, while it is contained in ${\rm FS}(\tilde f^-)$. Hence, 
$\alpha_{22}^l=0$. Analogously one obtains
\begin{equation}
\label{eq:f}
\alpha_{11}^r=\alpha_{22}^l=\beta_{12}^r=\beta_{21}^l=0.
\end{equation}
\begin{figure}[ht]
\begin{center}
\epsfysize=7cm
\epsfbox{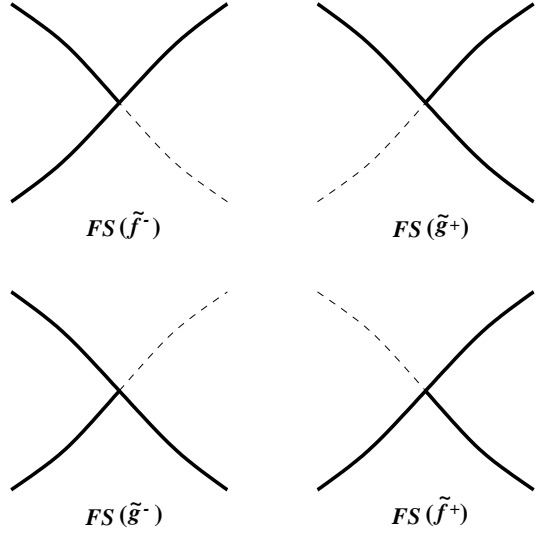}
\end{center}
\caption{Frequency sets of the microlocal solutions $\tilde f^\pm$ and $\tilde g^\pm$
defined in \eq{eq:tilde}.}
\end{figure}

%%%%%%%%%%%%%%%%%%%%%%%%%%%%%%%%%%%%%%%%%%%%%%%%%%%%%%%%%%%%%%%%%%%%%%%%%%%

\subsection{Connection between microlocal and exact WKB solutions}
\label{sec:con}
We now compute the remaining coefficients of the matrices $A_{l,r}$ and $B_{l,r}$ 
connecting the microlocal solutions $\tilde f^\pm$, $\tilde g^\pm$ with the WKB solutions
$u^\pm_l$, $u^\pm_r$ of $P_\nu u=E u$, $E>0$.

\medskip
As a first step, the stationary phase method gives  
the following formulae for the microlocal solutions, whose proof is to be found in Appendix \ref{app:stata}.

\begin{lem}\label{stata}
For the microlocal solutions $\tilde f^\pm$ of $P_\nu u=Eu$, $E>0$,  defined in \eq{eq:tilde}, one has microlocally near the Lagrangian manifolds $\sigma_l^+$ and $\sigma_r^-$ defined in \eq{eq:lag}
\begin{eqnarray}
\label{f-}
\tilde f^-(x) &=& 
e^{i\pi/8} \,2^{1/4} e^{+z(x;r_1)/h}
\left(\begin{array}{l} 0\\ 1\end{array}\right)
\left(1+O(h)\right)
\quad\mbox{near}\;\sigma_l^+,\\
\label{f+}
\tilde f^+(x)&=&
e^{i\pi/8}\,2^{1/4}  e^{-z(x;r_2)/h}
\left(\begin{array}{l} 0\\ 1\end{array}\right)
\left(1+O(h)\right)
\quad\mbox{near}\;\sigma_r^-.
\end{eqnarray}
\end{lem}

It remains to connect $e^{\pm z(x;r_1)/h}$ and $e^{\pm z(x;r_2)/h}$ to the exact WKB solutions $u_l^\pm$ and $u_r^\pm$.

\begin{prop} With the notation of Lemma \ref{stata}, we have microlocally
\label{relation}
$$
u^+_l=k^+_l\tilde f^-\quad\mbox{near}\,\,\,\sigma_l^+,
\quad
u^+_r=k^+_r\tilde g^+\quad\mbox{near}\,\,\,\sigma_r^+,
$$
$$
u^-_l=k^-_l\tilde g^-\quad\mbox{near}\,\,\,\sigma_l^-,
\quad
u^-_r=k^-_r\tilde f^+\quad\mbox{near}\,\,\,\sigma_r^-,
$$
where
$$
k_l^+=-2^{3/4}e^{i\pi/8}\left(1+O(h)\right),\qquad
k_r^-=-2^{3/4}e^{-i3\pi/8}\left(1+O(h)\right),
$$
and 
$$
k_l^-=-i\,\overline{k_l^+}
,\qquad
k_r^+=i\,\overline{k_r^-}
.
$$
\end{prop}

\begin{pr}
To prove these relations, say the first one, it is enough to calculate the asymptotic
behavior of $\tilde f^-$ and $u^+_l$ near $\sigma_l^+$.
First, recall the exact WKB formula 
$$
u^\pm_l(x)=e^{\pm z(x;r_1)/h} \,T_\pm(z(x;r_1)) 
\begin{pmatrix}w^{\rm even}_\pm(x)\\ w^{\rm odd}_\pm(x) \end{pmatrix}=
e^{\pm z(x;r_1)/h} \,T_\pm(z(x;r_1)) 
\begin{pmatrix}1\\ 0 \end{pmatrix}\big(1+O(h)\big)
$$
with
$$
T_\pm(z)=
\begin{pmatrix}
H^{-1}(z)\mp i H(z) & H^{-1}(z)\pm i H(z)\\
-H^{-1}(z)\mp i H(z) & -H^{-1}(z)\pm i H(z)
\end{pmatrix}
$$
and
$$
H(z(x))=\left(\frac{\nu+Ex-x^3}{\nu-Ex+x^3}\right)^{1/4}.
$$

\medskip
Let $x\in ]r_0, r_1[$ be fixed. With the branches of the fourth root chosen in Section \ref{global}, one has
$H(z(x))=e^{-i\pi/4}+O(h)$ as $h\to 0$. Hence, 
\begin{eqnarray*}
T_+(z(x))\left (
\begin{array}{c}
1 \\
0
\end{array}
\right )&=&\left (
\begin{array}{c}
H^{-1}(z(x))-iH(z(x)) \\
-H^{-1}(z(x))-iH(z(x))
\end{array}
\right )= -2\,e^{i\pi/4}\left (
\begin{array}{c}
0 \\
1
\end{array}
\right )
\left(1+O(h)\right)
,\\
T_-(z(x))\left (
\begin{array}{c}
1 \\
0
\end{array}
\right )&=&\left (
\begin{array}{c}
H^{-1}(z(x))+iH(z(x)) \\
-H^{-1}(z(x))+iH(z(x))
\end{array}
\right )= +2\,e^{i\pi/4}\left (
\begin{array}{c}
1 \\
0
\end{array}
\right )
\left(1+O(h)\right),
\end{eqnarray*}
and therefore
\begin{eqnarray}
\label{l+}
u_l^+(x)&=& -2\,e^{i\pi/4}e^{+z(x;r_1)/h}
\left(\begin{array}{c} 0\\1 \end{array}\right)
\left(1+O(h)\right),\\
u_l^-(x)&=& +2e^{i\pi/4}e^{-z(x;r_1)/h}
\left(\begin{array}{c} 1\\0 \end{array}\right)
\left(1+O(h)\right).
\nonumber
\end{eqnarray}
Comparing \eq{l+} and \eq{f-}, we immediately obtain 
$$
k_l^+=-2^{3/4}e^{i\pi/8}\left(1+O(h)\right).
$$

\medskip
Similarly, for fixed $x$ in the interval $]r_2,\infty[$, $H(z(x))=e^{\pi
i/4}+O(h)$ as $h\to 0$, and 
\begin{eqnarray}
\nonumber
u_r^+(x)&=& +2e^{-i\pi/4}e^{+z(x;r_2)/h}
\left(\begin{array}{c} 1\\0 \end{array}\right)
\left(1+O(h)\right),\\
\label{r-}
u_r^-(x)&=& -2e^{-i\pi/4}e^{-z(x;r_2)/h}
\left(\begin{array}{c} 0\\1 \end{array}\right)
\left(1+O(h)\right).
\end{eqnarray}
Comparing \eq{r-} and \eq{f+}, we immediately get 
$$
k_r^-=-2^{3/4}e^{-i3\pi/8}\left(1+O(h)\right).
$$

\medskip
For the computation of $k_r^+$ and $k_l^-$ we use some
symmetry properties. 
Recall that 
$$
g^{\pm}(y) =\left(\begin{array}{cc} 0&1\\1&0\end{array}\right)
{\mathcal C}{\mathcal
F}_hf^{\pm}(y), 
$$
and 
$$
{\mathcal C}{\mathcal F}_h ={\mathcal F}^{-1}_h{\mathcal C},
\quad{\mathcal C} V^{-1}=V {\mathcal C},
\quad V^{-1}{\mathcal F}^{-1}_h=V,
$$
where the last identity comes from ${\mathcal F}_h^{-1}$ being the metaplectic operator of the transformation $(x,\xi)\mapsto (\xi,-x)$, up to a normalizing constant.
Hence,
$$
V^{-1} g^{\pm}=
V^{-1}
\left(\begin{array}{cc}0&1\\1&0\end{array}\right) 
{\mathcal C}  {\mathcal F}_h f^{\pm}=
\left(\begin{array}{cc} 0&1\\1&0\end{array}\right) 
V^{-1} {\mathcal F}_h^{-1}{\mathcal C} f^{\pm}=
\left(\begin{array}{cc}0&1\\1&0\end{array}\right)
{\mathcal C}  V^{-1} f^{\pm}
$$
and
\begin{equation}
\label{symmetry2}
\tilde g^\pm=\left(\begin{array}{cc}
0&1\\1&0\end{array}\right)\overline{\tilde f^\pm}
.
\end{equation}
On the other hand,
$$
u_r^+=i\left(\begin{array}{cc}
0&1\\1&0\end{array}\right)\overline{u_{r}^-},\qquad
u_l^-=-i\left(\begin{array}{cc}
0&1\\1&0\end{array}\right)\overline{u_{l}^+},
$$
which together with \eq{symmetry2} yields
$$
k_l^-=-i\,\overline{k_l^+}
,\qquad
k_r^+=i\,\overline{k_r^-}
.
$$
\end{pr}

Proposition \ref{relation} yields for the connection matrices $A_{l,r}$ and $B_{l,r}$ defined in \eq{eq:mat} 
$$
\alpha_{21}^l=k_l^+,\quad
\alpha_{12}^r=k_r^-,\quad
\beta_{22}^l=k_l^-,\quad
\beta_{11}^r=k_r^+.
$$
Combining these with the knowledge on the vanishing matrix elements \eq{eq:f} and the relations $A_l=RB_l$, $A_r=RB_r$, one obtains
\begin{equation}
\label{eq:alr}
B_l=\begin{pmatrix}-\frac{1}{q}k_l^+&-\frac{p}{q} k_l^-\\*[1ex] 
0&k_l^-\end{pmatrix},\qquad
B_r=\begin{pmatrix}k_r^+&0\\*[1ex] -\frac{p}{q}k_r^+&\frac{1}{q}k_r^-\end{pmatrix}.
\end{equation}
%%%%%%%%%%%%%%%%%%%%%%%%%%%%%%%%%%%%%%%%%%%%%%%%%%%%%%%%%%%%%%%%%%%%
%%%%%%%%%%%%%%%%%%%%%%%%%%%%%%%%%%%%%%%%%%%%%%%%%%%%%%%%%%%%%%%%%%%%

\subsection{Computation of $T_2(E,h)$}
\label{sec:comp}

Finally, we prove Proposition \ref{7} calculating the  asymptotic behavior of $t(E,h)$ and $s(E,h)$ as $h\to0$. 

\medskip
Let $E>0$. By \eq{AB} and \eq{eq:alr},
\begin{equation}
\label{eq:t2}
T_2=B_r^{-1} B_l=
\begin{pmatrix}
{\displaystyle-\frac1{q}\frac{k_l^+}{k_r^+}}&
-{\displaystyle\frac{p}{q}\frac{k_l^-}{k_r^+}}\\*[3ex]
-{\displaystyle\frac{p}{q}\frac{k_l^+}{k_r^-}}&
{\displaystyle\frac{q^2-p^2}{q}\frac{k_l^-}{k_r^-}}
\end{pmatrix}.
\end{equation}

\begin{rem}
\label{check}
The identity \eq{eq:t2} is consistent with $T_2=\begin{pmatrix}t&s\\-\overline s&-\overline t\end{pmatrix}$,
since
$$
\overline{\left(\frac{k_l^+}{k_r^+}\right )}=-\frac{k_l^-}{k_r^-},\qquad
\overline{\left(\frac{k_l^+}{k_r^-}\right )}=-\frac{k_l^-}{k_r^+},\qquad
\overline{\left(\frac pq\right )}=\frac pq,\qquad
\frac 1{\overline q}=\frac{p^2-q^2}q,
$$
which is checked by direct calculation.
\end{rem}

By Proposition \ref{branching} and Proposition \ref{relation}, we then get
$$
t=-\frac 1q\frac{k_l^+}{k_r^+}=
-\frac{\sqrt\pi\,\gamma\,
\exp\!
\left(\frac{\pi}{4h}|\gamma|^2+\frac{i}{2h}|\gamma|^2\ln h\right)}
{\sqrt h\; \Gamma(1-\frac{i}{2h}|\gamma|^2)\;e^{i\pi/4}}
\left(1+O(h)\right).
$$
Since $\gamma=\frac{\tilde\nu}{\sqrt2}E^{-3/4}h+O(h^2)$ and
$$
\Gamma(1-{\textstyle\frac{i}{2h}}|\gamma|^2)=\sqrt{
\frac{\frac\pi{2h}|\gamma|^2}{\sin(\frac\pi{2h}|\gamma|^2)}
\left(1+{\textstyle\frac{|\gamma|^2}{4h^2}}\right)}=1+O(h),
$$ 
we have
$$
t=-\sqrt{\frac{\pi h}2}\,\tilde\nu\, E^{-3/4}\,e^{-i\pi/4}
+O(h|\ln h|)
$$
and
$$
s=-\frac pq\frac{k_l^-}{k_r^+}=
-i\,\exp\!\left({\textstyle\frac{\pi}{2h}}|\gamma|^2\right)
\left(1+O(h)\right)=
-i+O(h).
$$

\medskip
Since all the terms involved depend analytically on $E$ for $E\in\C$ near $E_0>0$, see Remark \ref{rem:gamma}, we have proven Proposition \ref{7}.

%%%%%%%%%%%%%%%%%%%%%%%%%%%%%%%%%%%%%%%%%%%%%%%%%%%%%%%%%%%%%%%%%%%%%%%%%%%%%
\section{Proof of the main results}
\label{distribution}

In this section, we compute the Bohr-Sommerfeld type quantization condition of Theorem \ref{main1}
and derive the semiclassical distribution of resonances given in Theorem \ref{main2}.

\subsection{Quantization condition}
Recall that Proposition \ref{quantizationcond} gives the quantization condition of resonances as $c^+(E,h)=0$, where 
$c^\pm(E,h)$ is the product of three transfer matrices
$T_1$, $T_2$, $T_3$ and the connection coefficients $c_0^\pm$,
$$
\left(\begin{array}{l}
c^+ \\
c^-
\end{array}
\right )=T_3T_2T_1
\left (
\begin{array}{l}
c^+_0 \\
c^-_0
\end{array}
\right ),
$$
see identity \eq{product}.
On the other hand, we have calculated the following asymptotics:
\begin{eqnarray*}
T_1&=&
\left (
\begin{array}{cc}
e^{S_{01}/h} & 0 \\
0 & e^{-S_{01}/h}
\end{array}
\right ),\\
T_2&=&\left (
\begin{array}{cc}
t & -i+O(h)\\
-i+O(h) & -\overline t
\end{array}
\right ),\quad
t=-\sqrt{\frac{\pi h} 2}e^{-i\pi/4}\tilde\nu E^{-3/4}+O(h|\ln h|),\\
T_3&=&2e^{-i\pi/4}
\left (
\begin{array}{cc}
e^{S_{2\infty}/h}\left(1+O(h)\right)  & O(e^{-\delta/h}) \\
O(e^{-\delta/h}) & e^{-S_{2\infty}/h}\left(1+O(h)\right) 
\end{array}
\right ),\\
\left (
\begin{array}{l}
c^+_0 \\
c^-_0
\end{array}
\right )&=&a
\left (
\begin{array}{l}
1+o(1) \\
-i+o(1)
\end{array}
\right ),
\end{eqnarray*}
see \eq{T1}, Proposition \ref{7}, Proposition \ref{prop:T3}, and Proposition \ref{airy}, respectively. 
Then,
$$
c^+=2a\,e^{-\pi i/4}\,e^{S_{2\infty}/h}\left(1+O(h)\right)
\left(t\,e^{S_{01}/h}(1+o(1))-e^{-S_{01}/h}(1+o(1))\right)
+O(e^{-\delta/h}).
$$
Hence, $c^+(E,h)=0$ if and only if
\begin{equation}
\label{eq:bsc}
\sqrt{\frac{\pi h}2}\,\tilde\nu\, e^{-i\pi/4}\, E^{-3/4}\, e^{2S_{01}(E,h)/h}+1=o(1)
\end{equation}
as $h\to0$, which proves Theorem \ref{main1}.

\subsection{Distribution of resonances}
We now study the asymptotic behavior of the function $S_{01}(E,h)$ as $h\to 0$. Recall that
$$
S_{01}(E,h)=\int_{r_0}^{r_1}\frac{\sqrt{\nu^2-r^2(E-r^2)^2}}{r} \,dr.
$$
If $\nu=\tilde\nu h>0$ is sufficiently small and $E>0$, then $\nu^2-r^2(E-r^2)^2>0$ for $r\in(r_0,r_1)$, and the square root in the formula for $S_{01}(E,h)$ is taken to be positive.
Substituting $y=r^2/E$, one gets
$$
S_{01}(E,h)=iE^{3/2}\int_{y_0}^{y_1}\frac{\sqrt{y(1-y)^2-\nu^2/E^3}}{2y} \,dy
$$
with $y_0=r_0^2/E$ and $y_1=r_1^2/E$.

\medskip
$y_0$ and $y_1$ are zeros of the cubic polynomial $y(1-y)^2-\mu^2$ with $\mu=\nu E^{-3/2}$. 
If $\mu>0$ is small and positive, then $y(1-y)^2-\mu^2$ has three  zeros $0<y_0(\mu)<y_1(\mu)<1<y_2(\mu)$ with $y_0(\mu)\to0$ and $y_{1,2}(\mu)\to1$ as $\mu\to0$.
We define 
$$
I(\mu)=\int_{y_0(\mu)}^{y_1(\mu)}\frac{\sqrt{y(1-y)^2-\mu^2}}{2y}\,dy,
$$
where the square root is taken to be positive for $0<\mu\ll1$.
Since
$$
S_{01}(E,h)=iE^{3/2}I(\mu),\qquad 
\mu=\frac{\nu}{E^{3/2}}=\frac{\tilde\nu}{E^{3/2}}h\,,
$$
we study the asymptotic behavior of the function $I(\mu)$
as $\mu\to 0$. For this, we have to understand the $\mu$-dependance of 
$y_0(\mu)$ and $y_1(\mu)$. When $\mu^2$ turns around $0$ once in the
positive sense (i.e. $\mu$ becomes $e^{\pi i}\mu$), then $y_0(\mu)$ turns around $0$ in the positive sense and $y_1(\mu)$ and $y_2(\mu)$ exchange their position turning half around $1$ in the positive sense. As a consequence, taking the branch into account,
\begin{equation}
\label{rot}
I(e^{\pi i}\mu)=I(\mu)+R(\mu)+T(\mu)
\end{equation}
with
$$
R(\mu)=-i\int_{\Gamma_0}\frac{\sqrt{\mu^2-y(1-y)^2}}{2y}dy,\quad
T(\mu)=-i\int_{y_1(\mu)}^{y_2(\mu)}\frac{\sqrt{\mu^2-y(1-y)^2}}{2y}dy,
$$
where $\Gamma_0$ is a contour around $0$.
These functions have the following properties:

\begin{lem}
\label{lem:tr}
$T(\mu)$ is a holomorphic function of $\mu^2$ at $\mu=0$, and
$$
R(\mu)=-\pi\mu,\quad
T(\mu)=\frac{\pi i\mu^2}{4}(1+O(\mu^2)).
$$
In particular, $R(e^{\pi i}\mu)=-R(\mu)$ and $T(e^{\pi i}\mu)=T(\mu)$. 
\end{lem}

\begin{pr}
By the residue theorem, $R(\mu)=\pi\sqrt{\mu^2}=-\pi\mu$.
For the study of $T(\mu)$, we move by the locally biholomorphic change of variables $v=\sqrt y(y-1)$ from a neighborhood of $y=1$ to a neighborhood of $v=0$,
$$
T(\mu)=-i\int_{-\mu}^{\mu}\sqrt{\mu^2-v^2}\,f(v)\,dv
$$
where $f(v)=(2y(v))^{-1}\frac{d}{dv}y(v)$ is holomorphic in a neighborhood of $v=0$ and satisifies 
$f(v)=\frac 12 +f'(0)v+O(v^2)$.
Hence,
$$
T(\mu)=i\mu^2\int_{-1}^1
\sqrt{1-w^2}\,f(\mu w)\,dw=
i\mu^2\left( 
{\textstyle\frac12}\int_{-1}^1 \sqrt{1-w^2}\,dw+O(\mu^2)
\right),
$$ 
since $\int_{-1}^1 \sqrt{1-w^2}\, w\, dw=0$.
\end{pr}

\begin{prop}
\label{las} 
$I(\mu)$ is ramified at $\mu=0$ and satisfies 
$$
I(\mu)={\textstyle\frac23}+{\textstyle\frac \pi 2} \mu+
O(\mu^2|\ln\mu|)
\qquad(\mu\to 0).
$$
\end{prop}

\begin{pr}
By Lemma \ref{lem:tr}, we have from (\ref{rot})
\begin{eqnarray*}
I(e^{2\pi i}\mu)&=&
I(e^{\pi i}\mu)+R(e^{\pi i}\mu)+T(e^{\pi i}\mu) =\left(I(\mu)+R(\mu)+T(\mu)\right)-R(\mu)+T(\mu) \\
&=&I(\mu)+2T(\mu).
\end{eqnarray*}
Since $\ln(e^{2\pi i}\mu)=\ln\mu+2\pi i$, 
this means that the function 
$$
A(\mu)=I(\mu)-{\textstyle\frac 1{\pi i}}\,T(\mu)\ln\mu
$$ 
is single-valued around $\mu=0$.
Moreover, since $T(\mu)$ behaves quadratically in $\mu$ near $\mu=0$, $A(\mu)$ is holomorphic near $\mu=0$ with
$$
A(\mu)\;\stackrel{\mu\to0}{\longrightarrow}\; 
A(0)=I(0)=\int_0^1(1-x^2)dx={\textstyle\frac23}.
$$ 
Differentiating equation (\ref{rot}), one gets
$A'(0)=I'(0)=\frac \pi 2$. Hence, 
\begin{eqnarray*}
I(\mu)&=&
A(0)+A'(0)\mu+{\textstyle\frac1{\pi i}}\,T(\mu)\ln\mu+O(\mu^2)=
{\textstyle\frac23}+{\textstyle\frac 1{\pi i}}\,T(\mu)\ln\mu+
{\textstyle\frac \pi 2}\mu
+O(\mu^2)\\
&=&{\textstyle\frac23}+{\textstyle\frac \pi 2} \mu+
O(\mu^2|\ln\mu|).
\end{eqnarray*}
\end{pr}

\begin{pr3}
The quantization condition \eq{eq:bsc} is satisfied, if and only if
there exists an integer $k\in\Z$ such that
\begin{equation}
\label{quant2}
2S_{01}(E,h)-{\textstyle\frac 12}\,h\ln E^{3/2}+
{\textstyle\frac 12}\,h\ln
\left ({\textstyle\frac{\pi}2}\tilde\nu^2h\right )=
\left (2k+{\textstyle\frac 54}\right )i\pi
h+o(h).
\end{equation}
Setting $\lambda=E^{3/2}$, Proposition \ref{las} implies
$$
S_{01}(E,h)=i\lambda\, I(\tilde\nu h/\lambda)=
{\textstyle\frac 2 3}\,i\lambda +
{\textstyle\frac {\pi} 2}\,i\tilde\nu h
+O(h^2|\ln h|),
$$
and the quantization condition
\eq{quant2} becomes
$$
{\textstyle\frac 4 3}i\lambda 
-{\textstyle\frac12}\,h\ln \lambda
+{\textstyle\frac 12}\,h\ln
\left ({\textstyle\frac{\pi}2}\tilde\nu^2 h\right )=\left (2k-\tilde\nu+{\textstyle\frac 54}\right
)i\pi h+o(h).
$$
Writing $\lambda=\lambda_1+i\lambda_2$ with $\lambda_1,\lambda_2\in\R$, the real and imaginary part of the above condition read as
\begin{eqnarray}
\label{quant4}
{\textstyle\frac 4 3}\lambda_1 
-{\textstyle\frac12}\,h\arg \lambda&=&\left (2k-\tilde\nu+{\textstyle\frac 54}\right
)\pi h+o(h),\\
\label{quant5}
-{\textstyle\frac 4 3}\lambda_2
-{\textstyle\frac12}\,h\ln |\lambda|+{\textstyle\frac 12}\,h\ln
\left ({\textstyle\frac{\pi}2}\tilde\nu^2 h\right )&=&o(h).
\end{eqnarray}
Now, we assume that $a<\lambda_1<b$ and $\lambda_2=o(1)$ as $h\to 0$. Then,
$$
\arg \lambda=\arctan(\lambda_2/\lambda_1)=o(1),\quad
\ln|\lambda|=\ln\lambda_1+{\textstyle\frac
12}\ln(1+\lambda_2^2/\lambda_1^2)=\ln\lambda_1+o(1).
$$
Setting 
$\lambda_{k\tilde\nu}=\frac{3\pi}{16}(8k-4\tilde\nu+5)$, 
euqations \eq{quant4} and \eq{quant5} become
\begin{eqnarray*}
\lambda_1 
&=&{\textstyle\frac {3\pi}{16}}(8k-4\tilde\nu+5) h+o(h) = \lambda_{k\tilde\nu}h + o(h),\\
\lambda_2
&=&
-{\textstyle\frac 38}\left(
h\ln{\textstyle\frac1h}-
h\ln\frac{\pi\tilde\nu^2}{2\lambda_{k\tilde\nu}h}\right)+o(h).
\end{eqnarray*}
\end{pr3}

%%%%%%%%%%%%%%%%%%%%%%%%%%%%%%%%%%%%%%%%%%%%%%%%%%%%%%%%%%%%%%%%%%%%%%%%%%%%%
\appendix

\section{Energy surfaces}
\label{app:es}

Let $E\in\C$ and $\nu\in h(\Z+\frac{1}{2})$. The zeros of the function  
$\C^+\to\C$, $r\mapsto(E-r^2)^2-\nu^2/r^2$ are the roots of the sixth order polynomial
$r^6-2Er^4+E^2r^2-\nu^2$ in $r$, which lie in the right half-plane $\C^+=\{r\in\C;\,
\re(r)>0\}$. This polynomial has at most three different roots $r_0,r_1,r_2\in\C^+$ in the
right half-plane, whose squares are the roots $x_0,x_1,x_2\in\C$ of the cubic polynomial
$x^3-2Ex^2+E^2x-\nu^2$ in $x$. If $E\in\R$, then an easy criterion for real-valuedness of
the roots $x_0,x_1,x_2$ is the sign check of the polynomial discriminant
$$
D_3 = (x_0-x_1)^2\, (x_0-x_2)^2\,(x_1-x_2)^2 = 
\nu^2\,(4E^3-27\nu^2)
\,.
$$
The three roots are real if and only if $D_3\ge0$, that is iff $\nu^2\le 4E^3/27$. 
The roots are real and distinct, if and only if $\nu^2<4E^3/27$.

From Cardano's formula
\begin{eqnarray*}
x_0 &=& 2E/3 + S_+ + S_-,\\
x_{1,2} &=& 2E/3 -(S_+ + S_-)/2 \pm i\sqrt{3}\,(S_+ - S_-)/2
\end{eqnarray*}
with $S_\pm=\sqrt[3]{-E^3/27+\nu^2/4\pm\sqrt{-D_3/108}}$ one learns 
$x_0\to0$ and $x_{1,2}\to E$ as $h\to0$.

If one denotes by $\sqrt{E}\in\C^+$ the square root of $E\in\C$, which lies in the right 
half-plane, then 
$r_0\to0$ and $r_{1,2}\to\sqrt{E}$ as $h\to0$, see also Figure 4.

\begin{figure}[h]
\label{fig:es}

%\vspace*{-10em}
\includegraphics[width=0.75\textwidth]{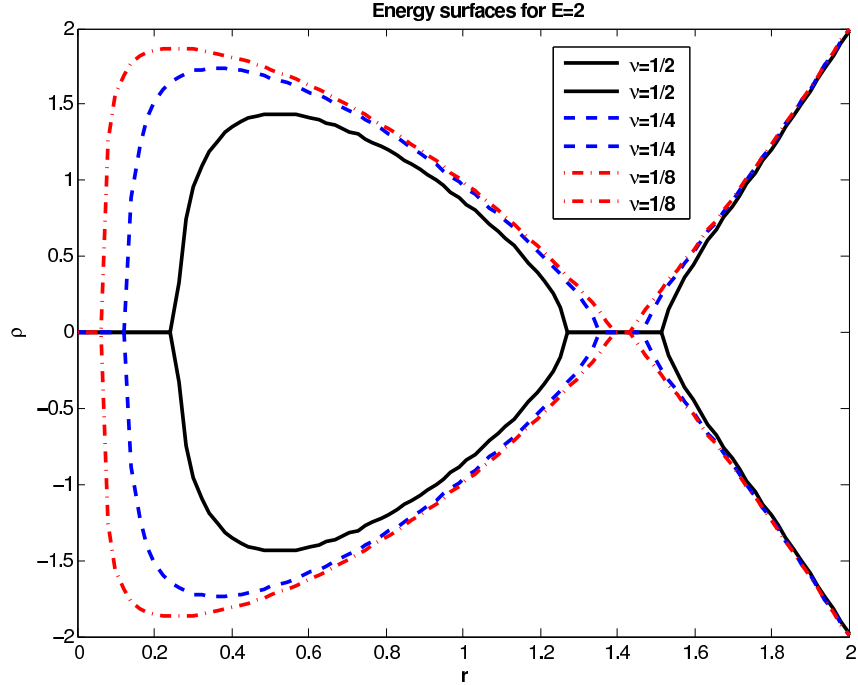}

%\vspace*{-10em}
\caption{The energy surfaces 
$\{(r,\rho)\in\R^+\times\R;\; \rho^2=(E-r^2)^2-\nu^2/r^2\}$ for $E=2$ and different values of 
$\nu\in\{\frac{1}{2},\frac{1}{4},\frac{1}{8}\}$.}
\end{figure}

%%%%%%%%%%%%%%%%%%%%%%%%%%%%%%%%%%%%%%%%%%%%%%%%%%%%%%%%%%%%%
\section{Spectrum of $P^+=-h^2\Delta + |x|$}
\label{app:plus}

The Schr\"odinger operator $P^+$ has a locally bounded positive potential, which increases 
to infinity as $|x|\to\infty$. Hence, $P^+$ is essentially self-adjoint on
$C^\infty_0(\R^2)$ and has purely discrete spectrum.  We are looking for eigenvalues
$E\in]a,b[$ in a bounded positive interval $]a,b[\subset\R^+$.  In polar coordinates
$x=r(\cos\theta,\sin\theta)$, $r>0$, $\theta\in\T$  the differential expression
$-h^2\Delta_x + |x|$ reads as
$$
-h^2\left(
\partial_r^2+\frac{1}{r}\partial_r+\frac{1}{r^2}\partial_\theta^2
\right)+r.
$$
Hence, $E$ is a {\em formal} solution of the eigenvalue problem $(P^+-E)\psi=0$ if and 
only if there exists $l\in\N\cup \{0\}$ such that 
$$
-h^2\left(
\frac{\d^2}{\d r^2}+\frac{1}{r}\frac{\d}{\d r}-\frac{l^2}{r^2}
\right)w_l(r)+(r-E)w_l(r)=0.
$$
Substituting $w_l(r)=r^{-1/2}\,u_l(r)$, this is equivalent to 
\begin{equation}
\label{eq:eig}
-h^2\left(
\frac{\d^2}{\d r^2}+\frac{\frac14-l^2}{r^2}
\right)u_l(r)+(r-E)u_l(r)=0.
\end{equation}

The ordinary differential equation \eq{eq:eig} has $r=0$ as a regular singular point 
with exponents $\frac12\pm l$, while $r=\infty$ is an irregular singular point of rank 
two. There are two linearly independent solutions, which behave as $r^{1/2}$ and 
$r^{1/2}\ln r$ (for $l=0$) or 
$r^{1/2+l}$ and $r^{1/2-l}$ (for $l\neq0$) near $r=0$, see for example Chapter 5.5 
in \cite{O}.
 
Looking at the irregular singularity at infinity, there is a fundamental system
$u^\pm_\infty(r)$, which behaves as 
$(r-E)^{-1/4} \exp(\pm(r-E)^{3/2}/h)$ near $r=\infty$. Hence, a  necessary condition for 
$E$ being an eigenvalue of $P^+$ reads as follows.

\begin{lem}
$E\in\sigma_{\rm disc}(P^+)$ if and only if  there exists $l\in\N\cup\{0\}$ such that $\W(u^0,u^-_\infty)=0$, 
where 
$\W(u^0,u^-_\infty)=u^0 \,(u^-_\infty)' - (u^0)'\, u^-_\infty$ 
is the Wronskian of the two solutions to \eq{eq:eig},
\begin{eqnarray*}
u^0(r)&\sim& r^{1/2+l},\qquad r\to0,\\
u^-_\infty(r)&\sim&(r-E)^{-1/4}\, \e^{-(r-E)^{3/2}/h},\qquad 
r\to\infty.
\end{eqnarray*}
\end{lem}
 
Let us have a closer look at the potential of equation (\ref{eq:eig})
$$
Q(r)=r-E+h^2(l^2-1/4)/r^2.
$$ 
The origin $r=0$ is a double pole and the coefficient $h^2(l^2-1/4)$ is positive for
$l>0$ and is negative for $l=0$. This difference is crucial from the WKB point of view, since 
the level curves
of $\re(\int\sqrt{Q(r)}\,dr)$ in the complex $r$-plane are closed curves enclosing the origin for $l>0$ and radial curves emanating the origin for $l=0$.
The WKB construction at a double pole for $l>0$ has been treated in \cite{langer},
and we here restrict ourselves to the application of this results. 

We assume $l>0$.
For $h>0$ sufficiently small, the potential $Q$
has three distinct simple turning points: they are zeros of the cubic polynomial 
$r^3-Er^2+h^2(l^2-1/4)$.   
Two of them are at a distance $O(h)$ from $r=0$, while the third is at a distance $O(h)$ 
from $r=E$.
Since the polynomial discriminant 
$$
h^2(l^2-{\textstyle\frac14})
\left(
{\textstyle\frac1{27}}\,E^3-
{\textstyle\frac{h^2}4}(l^2-{\textstyle\frac14})
\right)>0
$$
is positive, all three turning points are real.
The first one is negative, while the other two are positive.
We denote these two by $\alpha_1,\alpha_2>0$.
The strategy for characterizing the quantization condition 
$\W(u^0,u^-_\infty)=0$ is the following. We connect
\begin{enumerate}
\item
$u^0$ with two exact WKB solutions $\widetilde u^\pm$ built with the 
Langer-modified potential 
$$
\widetilde Q(r)=r-E+h^2 l^2/r^2
$$ 
defined  for $r\in]\alpha_1,\alpha_2[$ with the phase base point at $\alpha_1$,
\item
$\widetilde u^\pm$ with two exact WKB solutions $u^\pm_l$ defined for 
$r\in]\alpha_1,\alpha_2[$ with the phase base point at $\alpha_2$,
\item
$u^\pm_l$ with an exact WKB solution $u^-_r$ defined for $r>\alpha_2$, which is collinear to $u^-_\infty$.
\end{enumerate}

We start with the connection near the origin.
Proposition 11 in \cite{langer} proves existence of a non-zero constant $a=a(E,h)\neq0$ such that
$$
u^0 = 
a\left(1+o(1)\right)\widetilde u^+ + 
a\left(\i+o(1)\right)\widetilde u^-
$$
as $h\to0$, where 
$$
\widetilde u^\pm(r)=\widetilde Q(r)^{-1/4} 
\exp\!\left(
\pm\frac1h\int_{\widetilde\alpha_1}^r \widetilde Q(s)^{1/2}\d s
\right) 
\widetilde W^\pm(r;r^\pm_l)\,.
$$
$\widetilde\alpha_1>0$ denotes the first positive turning point of $\widetilde Q$, while $r^\pm_l\in\C$ are two suitably chosen $h$-independent points in the complex plane such that 
$\Re(r^\pm_l)\in]\alpha_1,\alpha_2[$, $\pm \Im(r^\pm_l)>0$, and $\widetilde W^\pm(r^\pm_l;r^\pm_l)=1$.

Next, we connect near the turning point.
The WKB solutions $u^\pm_l(r)=u^\pm(r;\alpha_2,r^\pm_l)$ are of the form
$$
u^\pm_l(r)=Q(r)^{-1/4} 
\exp\!\left(\pm\frac1h\int_{\alpha_2}^r Q(s)^{1/2}\d s\right) 
W^\pm(r;r^\pm_l)
$$
with $W^\pm(r^\pm_l;r^\pm_l)=1$. 
Connecting $\widetilde u^\pm$ and $u^\pm_l$, 
$$
(\widetilde u^+,\widetilde u^-)=(u^+_l,u^-_l)
\begin{pmatrix}
\W(\widetilde u^+,u^-_l) & \W(\widetilde u^-,u^-_l)\\ 
-\W(\widetilde u^+,u^+_l) & -\W(\widetilde u^-,u^+_l)
\end{pmatrix}
\W(u^+_l,u^-_l)^{-1},
$$
we evaluate the Wronskians, which are constant with respect to $r$, near the points $r=r^\pm_l$. 
By this choice we can exploit the fact that 
$\widetilde Q(r) = Q(r) + O(h^2)$  
as $h\to0$ uniformly in neighborhoods of $r=r^\pm_l$. Also the turning points differ by a term of order $h^2$, $\widetilde\alpha_1=\alpha_1+O(h^2)$, and one derives
$$
\widetilde u^\pm(r) = u^\pm(r;\alpha_1,r^\pm_l) \left(1+O(h)\right)
$$
as $h\to0$ uniformly in neighborhoods of $r=r^\pm_l$. Hence for the connection formula, we just have to account for solutions with different phase base points $\alpha_1$, $\alpha_2$.
By Proposition 2.5 in \cite{R}, one has
\begin{eqnarray*}
\W(u^+_l,u^-_l)&=&-{\textstyle\frac2h}\left(1+O(h)\right),\\
\W(\widetilde u^\pm,u^\pm_l)&=&
O(\e^{-\delta/h}),\\
\W(\widetilde u^\pm,u^\mp_l)&=&
\mp{\textstyle\frac2h}\,\e^{\pm S_{12}/h}\left(1+O(h)\right),
\end{eqnarray*}
where $\delta>0$ does not depend on $h$, and 
$$
S_{12}=S_{12}(E,h)=\int_{\alpha_1}^{\alpha_2}  
Q(r)^{1/2}\,\d r
$$
is the action integral between $\alpha_1$ and $\alpha_2$. We note, that the square root is taken such that $S_{12}$ is purely imaginary with positive imaginary part.
Collecting the previous formulae, one obtains
$$
(\widetilde u^+,\widetilde u^-)=(u^+_l,u^-_l)
\begin{pmatrix}
\e^{S_{12}/h}\left(1+O(h)\right) & O(\e^{-\delta/h})\\ 
O(\e^{-\delta/h}) & \e^{-S_{12}/h}\left(1+O(h)\right)
\end{pmatrix}
$$
and
\begin{equation}
\label{eq:con}
u^0 = a\,\e^{S_{12}/h}\left(1+o(1)\right) u^+_l+
a\,\e^{-S_{12}/h}\left(\i+o(1)\right) u^-_l\qquad (h\to0).
\end{equation}

Finally, we connect to infinity.
The solution $u^-_r(r)=u^-_r(r;\alpha_2,r_r)$ is defined as
$$
u^-_r(r)=Q(r)^{-1/4} 
\exp\!\left(-\frac1h\int_{\alpha_2}^r Q(s)^{1/2}\d s\right) 
W^-(r;r_r)
$$
with $r_r\in\R$ such that $r_r>\alpha_2$ and $W^-(r_r;r_r)=1$. 
Since 
$$
Q(r)\sim r-E,\qquad r\to\infty,
$$
$u^-_r$ is collinear to $u^-_\infty$, and $\W(u^0,u^-_\infty)=0$ if and only if $\W(u^0,u^-_r)=0$.  By equation (\ref{eq:con}), it is enough to compute the Wronskians 
$\W(u^\pm_l,u^-_r)$. 
Choosing a branch cut for $Q(r)^{1/2}$ near $\alpha_2$ as in Figure 2 of \cite{langer}, one gets
$$
\W(u^+_l,u^-_r)=-{\textstyle\frac2h}\left(1+O(h)\right),\qquad
\W(u^-_l,u^-_r)={\textstyle\frac{2\i}h}\left(1+O(h)\right),
$$
since $u^-_l(r)=-\i\, u^+(r;\alpha_2,s^-_l)$, where 
$s^-_l$ is the point $r^-_l$ on the other Riemann surface. 
Hence, we have proven

\begin{prop}
\label{prop:plus}
Assume $l>0$. Then, 
$$
{\mathcal W}(u^0,u^-_\infty)=0\quad\Longrightarrow\quad 
\e^{2S_{12}(E,h)/h}+1=o(1),\quad h\to0,
$$
with  
\begin{equation}
\label{eq:act}
S_{12}(E,h)=\int_{\alpha_1}^{\alpha_2} 
\sqrt{r-E+h^2(l^2-{\textstyle\frac14})/r^2}\;\d r.
\end{equation}
\end{prop}

What remains, is a study of the asymptotic behavior of the action integral $S_{12}(E,h)$ as $h\to0$. Substituting $r=Ey$ in \eq{eq:act},  
we rewrite 
$$
S_{12}(E,h)=\i\,E^{3/2}I^+(\mu),\qquad
\mu = \frac{\sqrt{l^2-1/4}}{E^{3/2}}\,h,
$$ 
with
$$
I^+(\mu)=\int_{y_1(\mu)}^{y_2(\mu)} 
\frac{\sqrt{y^2(1-y)-\mu^2}}{y}\,\d y,
$$
where $y_1(\mu)$ and $y_2(\mu)$ are the first and second positive zero of the cubic polynomial $y^2(1-y)-\mu^2$. For $\mu>0$ small enough, the polynomial has three zeros $y_0(\mu)<0<y_1(\mu)<y_2(\mu)<1$ with $y_0(\mu),y_1(\mu)\to0$ and $y_2(\mu)\to1$ as $\mu\to0$. When $\mu$ turns half around zero in the positive sense, that is when $\mu$ turns to $\e^{\i\pi}\mu$, then $y_0(\mu)$ and $y_1(\mu)$ exchange their position, while $y_2(\mu)$ turns around $1$. Hence,
$I^+(\e^{\i\pi}\mu) = I^+(\mu) + R^+(\mu) + T^+(\mu)$
with
$$
R^+(\mu)=-\i\int_{y_1(\mu)\curvearrowright y_0(\mu)}\;
\frac{\sqrt{\mu^2-y^2(1-y)}}{y}\;\d y
$$
and $T^+(\mu)=\int_{\Gamma_1}\sqrt{y^2(1-y)-\mu^2}/y\,\d y=0$,
where $\Gamma_1$ is a contour around $1$. Substituting $v=y\sqrt{1-y}$, we rewrite
$$
R^+(\mu)=-\i\int_{\mu\curvearrowright -\mu}
\sqrt{\mu^2-v^2}\;\frac{f(v)}v\;\d v,
$$ 
where $f(v)=v^2/(v^2-\frac12 y(v)^3)$ is holomorphic in a neighborhood of $v=0$ and satisfies $f(v)=1+\frac12 v + f''(0)v^2 + O(v^3)$.
Then,
$$
R^+(\mu)=-\i\,\mu\int_{1\curvearrowright -1}
\sqrt{1-w^2}\;\frac{f(\mu w)}w\;\d w=
-\i\,\mu \left(-\i\pi -{\textstyle\frac{\pi}{4}}\mu + O(\mu^3)\right)
$$
as $\mu\to0$, 
since $\int_{1\curvearrowright -1}
\sqrt{1-w^2} \,w^j \,\d w$ equals $-\i\pi,-\frac\pi2,0$ for 
$j=-1,0,1$, respectively. Hence,
$$
I^+(\e^{\i\pi}\mu) = I^+(\mu) -\pi\mu +{\textstyle\frac{\pi\i}{4}}\,\mu^2 \left(1+O(\mu^2)\right)\qquad(\mu\to0),
$$
and Proposition \ref{las} yields
$$
I^+(\mu)={\textstyle\frac23 +\frac\pi2}\mu + O(\mu^2|\ln\mu|)\qquad(\mu\to0).
$$

\begin{rem}
With these asymptotics, 
the quantization condition of Proposition \ref{prop:plus} can be rephrased as follows. There exist an integer $k\in\Z$ and a natural number $l\in\N$ such that
$S_{12}(E,h)=\i\pi (k+\frac12) h +o(h)$ as $h\to0$ or 
$$
E^{3/2}={\textstyle\frac34}\pi\left( 2k + 1 - \sqrt{l^2-1/4}\right) h + o(h)
\qquad (h\to0). 
$$
\end{rem}

The study of the resonant set of the full operator $P$ and the spectrum of the upper level operator $P^+$ have revealed a parallel structure: A resonance $E$ of $P$ is characterized by a Bohr-Sommerfeld type condition
$$
\tilde t\,\e^{2 S_{01}(E,h)/h} + 1 = o(1)\qquad(h\to0)\,,
$$ 
where $\tilde t$ is related to non-adiabatic transitions, and $S_{01}(E,h)$ is an action integral. An eigenvalue $\lambda$ of $P^+$ is characterized by a Bohr-Sommerfeld condition
$$
\e^{2 S_{12}(\lambda,h)/h} + 1 = o(1)\qquad(h\to0)\,,
$$
where $S_{12}(\lambda,h)$ is an action integral. The integrals can be expressed as
\begin{eqnarray*}
S_{01}(E,h) = \i\, E^{3/2} I(\mu),\qquad 
\mu = h\, \tilde\nu\, E^{-3/2},\qquad
\tilde\nu\in\N+{\textstyle\frac12},\\
S_{12}(\lambda,h) = \i\, \lambda^{3/2} I^+(\mu),\qquad 
\mu = h \,\sqrt{l^2-{\textstyle\frac14}} \,\lambda^{-3/2},\qquad
l\in\N,
\end{eqnarray*}
where $I(\mu)$ and $I^+(\mu)$ share the same $\mu$-asymptotics for $\mu\to0$.

%%%%%%%%%%%%%%%%%%%%%%%%%%%%%%%%%%%%%%%%%%%%%%%%%%%%%%%%%%%%
\section{Symbols, frequency sets, microlocal solutions}
\label{app:fs}

Let us briefly recall the definition of the symbol classes used, as well as the notion of frequency
set and microlocal solution.

\begin{defi}
A smooth function $a(x,h)$ on an open set $\Omega\subset\R$ is
called a 
$\coi$-{\em symbol} of order $m\in\Z$, if there exist smooth functions $a_n(x)\in \coi(\Omega,\C)$ such that
for all non-negative integers $\alpha\in\N$ and $N\in\N$,
$$
\label{symbol1}
\sup_{x\in\Omega}\mid\partial^{\alpha}_{x}(a(x,h)-\sum_{n=0}^{N}a_n(x)h^{n+m})\mid
\;=\;O(h^{m+N+1}),\qquad h\to0.
$$
When this holds, we write 
$a(x,h)\sim\sum_{n=0}^\infty a_n(x) h^{n+m}$.
\end{defi}

\begin{defi}
An analytic function $a(x,h)$ on an open set $\Omega\subset\C$ is called
a {\em Gevrey symbol} of index $2$, if there exist functions $a_n(x)$ analytic in $\Omega$ such that
for all compact subsets $K\subset\Omega$ and all large enough $C>0$ there exists $\delta>0$ with
$$
\sup_{x\in K}\mid a(x,h)-\sum_{n=0}^{\frac{1}{C\sqrt{h}}}a_n(x)h^{n}\mid
\;=\;O(e^{-\delta/\sqrt{h}}),\qquad h\to0.
$$   
\end{defi}

\begin{defi}
Let $u\in{\mathcal S}'$ be a possibly $h$-dependent distribution 
and $(x_0,\xi_0)\in T^*\R$.
Then $(x_0,\xi_0)\notin {\rm FS}(u)$ if there exists a Schwartz function $\chi_0$ on phase space $T^*\R$ with $\chi_0(x_0,\xi_0)=1$ such that for any $N\in\N$
$$
\chi_0(x,hD_x) u=O(h^N),\qquad h\to0.
$$
${\rm FS}(u)$ is called the {\em frequency set} of $u$, and 
$u$ a {\em microlocal solution of} $u=0$ near $(x_0,\xi_0)$. 
\end{defi}

In Sections \ref{coni} and \ref{sec:con}, we have used the following fact about the frequency set of WKB solutions:

\begin{lem}
\label{bkw}
If $u(x,h)=a(x,h)\exp (i\phi(x)/h)$ where $\phi\in {\mathcal C}^1(\Omega,\R)$ is a real phase and $a$ is a 
${\mathcal C}^\infty$-symbol, then
${\rm FS}(u)\subset \{(x,\xi);\; \xi=\partial_x\phi(x)\}$.
\end{lem}

%%%%%%%%%%%%%%%%%%%%%%%%%%%%%%%%%%%%%%%%%%%%%%%%%%%%%%%%%%%%%%%%%%%%%%%%%%
\section{Proof of Theorem \ref{tono}}
\label{app:pr}

\begin{lem}
\label{lem:l}
Let $\tilde\nu\in \N-\frac12$, $y\mapsto\psi(y)$ the function defined in \eq{eq:psi} with
$\psi(0)=\frac{E^{-3/4}}{\sqrt2}$, $E>0$, and $v(y)=v(y,h)$ a solution of
\begin{equation}
\label{eq:v}
h D_y v(y)=
\begin{pmatrix}y&h\,\tilde\nu\psi(y)\\
-h\,\tilde\nu\psi(y)&-y\end{pmatrix}
v(y).
\end{equation} 
There exists a matrix-valued ${\mathcal C}^\infty$-symbol $M(y,h)={\rm Id}+O(h)$, such that $w(y,h)=M(y,h)v(y,h)$ satisfies 
\begin{equation}
\label{eq:w}
\begin{pmatrix}
hD_y-y&-\gamma\\
\overline\gamma&hD_y+y
\end{pmatrix}
w(y,h)=r(y,h)w(y,h),
\end{equation}
where $\gamma=\frac{\tilde\nu}{\sqrt2}E^{-3/4} h+O(h^2)$ and $r(y,h)=O(h^\infty)$ uniformly in an interval around $y=0$ together with all its derivatives.
\end{lem}

\begin{pr}
We rewrite equation \eq{eq:v} as $hD_y v(y,h)=B(y,h)v(y,h)$ with 
$$
B(y,h)=B_0(y)+hB_1(y)=
\begin{pmatrix}y&0\\0&-y\end{pmatrix}+
h\begin{pmatrix}0 & \tilde\nu\psi(y) \\
-\tilde\nu\psi(y) & 0\end{pmatrix}
$$
and equation \eq{eq:w} as $\left(hD_y-G(y,h)\right)w(y,h)=r(y,h)w(y,h)$ with 
$$
G(y,h)\sim \sum_{n=0}^\infty G_n(y) h^n,\qquad
G_0(y)=B_0(y),\qquad
G_n(y)\equiv 
\begin{pmatrix}0&\gamma_n\\-\overline\gamma_n & 0\end{pmatrix}
\qquad (n\ge1).
$$
We are looking for 
$$
M(y,h)\sim \sum_{n=0}^\infty M_n(y)h^n,
\qquad
M_n(y)=
\begin{pmatrix}
m_n(y) & q_n(y) \\*[1ex]
\overline q_n(y) & \overline m_n(y)
\end{pmatrix}
\qquad(n\ge0)
$$
such that
$$
hD_yM=GM-MB
$$
or equivalently for all $n\ge 0$
\begin{equation}
\label{rec}
D_yM_{n-1}=\sum_{j=0}^n\left(G_jM_{n-j}-M_jB_{n-j}\right)
\end{equation}
with convention $M_{-1}=0$, $B_n=0$ for $n\ge 2$.

Equation \eq{rec} is satisfied for $n=0$, if we take 
$M_0(y)\equiv{\rm Id}$, i.\ e.\  $m_0(y)\equiv 1$ and $q_0(y)\equiv 0$. 

Then, we have for $n\ge 1$
\begin{eqnarray*}
G_n &=& -i M_{n-1}' - \sum_{j=0}^{n-1} G_j
M_{n-j} + M_{n-1}  B_1 + M_n  B_0\\
&=& -i M_{n-1}' - \sum_{j=1}^{n-1} G_j
M_{n-j} + M_{n-1}  B_1 -2y\left(\begin{array}{cc}0&
q_n\\ - \overline q_n&0\end{array}\right).
\end{eqnarray*}

Since $\psi(y)\in\R$ for $y\in\R$, the previous equation is equivalent to
\begin{eqnarray}
\label{indu1}
\gamma_n &=&- iq_{n-1}'(y) -
\sum_{j=1}^{n-1}\gamma_j \overline m_{n-j}(y) +
\tilde{\nu}\psi(y)m_{n-1}(y) - 2 yq_n(y),\\
\label{indu2}
0 &=& -im_{n-1}'(y) -
\sum_{j=1}^{n-1}\gamma_j \overline q_{n-j}(y) -
\tilde{\nu}\psi(y)q_{n-1}(y).
\end{eqnarray}
Let us start with $n=1$. 
Substituting $y=0$ in \eq{indu1} yields
$$
\gamma_1=\tilde\nu\psi(0),
$$
and we automatically obtain
$$
q_1(y)=-\frac 1{2y}\left(\gamma_1-\tilde\nu\psi(y)\right),
$$
which is smooth (even analytic) near $y=0$. Setting $m_1(0)=0$, equation \eq{indu2} gives
$$
m_1(y)=i\int_0^y \left(\gamma_1 \overline q_1(y') +\tilde\nu\psi(y') q_1(y')\right)dy'.
$$
If $\gamma_j,m_j,q_j$ with $m_j(0)=0$ are determined for $1\le j\le n$, then \eq{indu1} yields
$$
\gamma_{n+1}=-i q_n'(0)
$$
and
$$
q_{n+1}(y)=-\frac
1{2y}\big(
\gamma_n+iq_n'(y)+\sum_{j=1}^n\gamma_j\overline m_{n+1-j}(y)-\tilde
\nu\psi(y)m_n(y)
\big),
$$
which is smooth (even analytic) near $y=0$.
Setting $m_{n+1}(0)=0$, equation \eq{indu2} then gives
$$
m_{n+1}(y)=i\int_0^y\big(\sum_{j=1}^{n+1}\gamma_j\overline q_{n+2-j}(y')+
\tilde\nu\psi(y')q_{n+1}(y')\big)dy'.
$$

Since $m_n(y)$ and $q_n(y)$, $n\ge0$, are smooth functions near $y=0$, there exists a matrix-valued $\coi$-symbol $M(y,h)$ with $M(y,h)\sim \sum_{n\ge0}M_n(y) h^n$ such that \eq{eq:w} holds.
\end{pr}

\begin{lem}
\label{lem:meta}
Let 
$\kappa_{\frac\pi 4}(y,\eta)={\textstyle\frac1{\sqrt2}}(y-\eta,y+\eta)$ 
be the $\frac{\pi}4$-rotation in phase space $T^*\R$. The metaplectic operator $V$ of the transpose linear canonical transformation 
$\kappa_{\frac\pi 4}^*=\kappa_{-\frac\pi 4}$ satisfies
$$
V(h D_y -y) = -\sqrt2 y V\,,\qquad V(h D_y + y) = \sqrt2 \,h D_y V.
$$

\noindent
Moreover, ${\rm FS}(Vu)=\kappa_{-\frac\pi 4}{\rm FS}(u)$ for $u\in{\mathcal S}'$.
\end{lem}

\begin{pr}
One applies Theorem 2.15 in \cite{fo} for the linear canonical transformation $\kappa_{\frac\pi 4}$.
\end{pr} 

\begin{rem}
A formula for $V$ in terms of oscillatory integrals is given by
$$
V g(y) = e^{i\pi/8}(\sqrt2\pi h)^{-1/2}\int_\R e^{-\frac i{2h} (y^2-2\sqrt 2 x y + x^2)} g(x) dx, 
$$
see Theorem 4.53 in \cite{fo} or Proposition 5.3 in \cite{R}.
\end{rem} 
%%%%%%%%%%%%%%%%%%%%%%%%%%%%%%%%%%%%%%%%

%%%%%%%%%%%%%%%%%%%%%%%%%%%%%%%%%%%%%%%
\section{Gevrey symbols}
\label{app:gev}

\begin{lem}
The analytic functions $m_0(x)\equiv1$, $q_0(x)\equiv0$, 
\begin{eqnarray*}
q_n(x) &=& -{1\over 2x} 
\Big(i\int_0^x q_{n-1}''(t)dt  
-i\sum_{j=2}^{n-1}q_{j-1}'(0)\,\overline m_{n-j}(x)\\
&&\hspace*{6em}
-m_{n-1}(x)\,\tilde \nu\,\psi(x)+\overline m_{n-1}(x)\,\tilde\nu\psi(0)\Big),
\qquad n\ge1,\\
m_n(x) &=& i\,\tilde{\nu} \int_0^x q_n(t) \psi(t) dt +
i\sum_{j=1}^n \gamma_{j}\int_0^x \overline q_{n+1-j}(t) dt,\qquad n\ge1\\
\end{eqnarray*}
with $\gamma_1=\tilde\nu \psi(0)$ and $\gamma_n= -i q_{n-1}'(0)$ for $n\ge 2$, which have been introduced in the proof of Lemma \ref{lem:l}, have analytic resummations $q(x,h)$ and $m(x,h)$ in some neighborhood of $x=0$, which are Gevrey symbols of index $2$.
\end{lem}

\begin{pr}
Let
$\Omega=\{x\in\C;\,|x|<r\}$ and $\Omega_t=\{x\in\C;\,|x|\le r-t\}$
for $0<t\le r$. It is enough to show, that 
there exist positive constants $D_m,D_q,C>0$ such that for all $n\in\N$
\begin{equation}
\label{pun}
\sup_{x\in \Omega_t}\left|q_n(x) \right| \leq 
D_q C^n \frac {{(2n)}^{2n}}{t^{2n}},\qquad
\sup_{x\in \Omega_t}\left|m_n(x)\right| \leq D_m C^n \frac
{{(2n)}^{2n}}{t^{2n}}.
\end{equation}

We use the following bounds: if $f(x)$ is a holomorphic function in $\Omega$, which satisfies for some positive constant $M>0$
$$
\sup_{x\in\Omega_t}|f(x)|\le \frac M{t^k},
$$
then 
$$
\sup_{x\in\Omega_t}\left |\int_0^xf(u)du\right |\le \frac M{(k-1)t^{k-1}},
\qquad
\sup_{x\in\Omega_t}\left |\frac 1x\int_0^xf(u)du\right |\le \frac M{t^{k}},
$$
and
$$
\sup_{x\in\Omega_t}|f'(x)|,\;
\sup_{x\in\Omega_t}\left|{\textstyle\frac{f(x)}x}\right|
\le \frac {(k+1)^{k+1}M}{k^k\, t^{k+1}},
\qquad
\sup_{x\in\Omega_t}|f''(x)|\le \frac {(k+2)^{k+2}M}{k^k\, t^{k+2}}.
$$

\bigskip
We prove \eq{pun} by induction:
\begin{eqnarray*}
\sup_{x\in \Omega_t}\left|q_n(x)\right|&\le&
\sup_{x\in \Omega_t}\left|q_{n-1}''(x)\right|+
\sum_{j=2}^{n-1}|q_{j-1}'(0)|\sup_{x\in
 \Omega_t}\left |{\textstyle\frac {m_{n-j}(x)}{x}}\right|\\
&& \hspace*{10em}
+ 2|\tilde\nu|\sup_{x\in \Omega_t}\left|{\textstyle\frac {m_{n-1}(x)} x} \right|\sup_{x\in\Omega}|\psi(x)|\\
&\le& 
D_q C^{n-1} \frac {{(2n)}^{2n}}{t^{2n}}+ 
\sum_{j=2}^{n-1}
D_q C^{j-1} \frac {{(2j-1)}^{2j-1}}{t^{2j-1}}
D_m C^{n-j} \frac {{(2n-2j+1)}^{2n-2j+1}}{t^{2n-2j+1}}\\
&& \hspace*{10em}
+2|\tilde\nu|D_m C^{n-1} \frac {{(2n-1)}^{2n-1}}{t^{2n-1}}
\sup_{x\in \Omega}|\psi(x)|\\
&\le&
D_q\left(1+D_m+2|\tilde \nu|\,r\,\frac{D_m}{D_q}\,\sup_{x\in \Omega}|\psi(x)|\right)
C^{n-1}
\frac {{(2n)}^{2n}}{t^{2n}},
\end{eqnarray*}
where we used the inequality
$$
k^k(N-k)^{N-k}\le \left ({\textstyle\frac N2}\right )^N\qquad (0<k<N).
$$
Taking $C$ large enough so that $C\ge 1+D_m +2 r|\tilde\nu| \frac
{D_m}{D_q}
\sup_{x\in \Omega}|\psi(x)|$, we obtain the first bound in \eq{pun}.
Similarly, one has
\begin{eqnarray*}
\sup_{x\in \Omega_t} \left |m_n(x)\right | &\le&
|\tilde{\nu}|\sup_{x\in \Omega}|\psi(x)|
D_q C^n  \frac {(2n)^{2n}} {(2n-1)t^{2n-1}}\\
&&
+\sum_{j=1}^n D_q C^{j-1}  {(2j-1)}^{2j-1}
D_q C^{n+1-j} \frac{(2n+2-2j)^{2n+2-2j}} {(2n+1-2j)\,t^{2n}}\\
&\leq& \big(|\tilde{\nu}|\,r\,\sup_{x\in \Omega}|\psi(x)| +2D_q \big)D_q 
C^n \frac {{(2n)}^{2n}}{t^{2n}}.
\end{eqnarray*}
Taking $D_m$ so that $D_m\ge (|\tilde{\nu}|r\sup_{x\in
\Omega}|\psi(x)|+  2D_q)D_q $, we get the second bound in \eq{pun}. 
\end{pr}

%%%%%%%%%%%%%%%%%%%%%%%%%%%%%%%%

\section{Proof of Lemma \ref{stata}}
\label{app:stata}

In this appendix,  we prove one of the two formulae of Lemma \ref{stata},
i.e. compute the asymptotic exapnsion of $\tilde f^+$ microlocally near
$\sigma_r^-$. 

\medskip
Recall that near $y=0$
$$
\tilde f^+(\phi^{-1}(y))=M(y,h)^{-1}V^{-1}f^+(y)
$$
where the inverse of the metaplectic operator is given by
$$
V^{-1}f^+(y)=e^{i\pi/8}(\sqrt 2\pi h)^{-1/2}
\int_\R e^{i(y^2+z^2-2\sqrt2yz)/(2h)}f^+(z)dz,
$$
and that
$$
f^+= {}^t(f_1^+,f_2^+), 
\quad
f_1^+= -\frac {\gamma}{\sqrt{2}y}\chi_{(0,\infty)}(
y)\,|y|^{\frac i {2h} |\gamma|^2 }\,,
\quad
f_2^+= \chi_{(0,\infty)}(y)\,|y|^{\frac i {2h} |\gamma|^2
}\,.
$$ 
Let us compute the asymptotic expansion of 
$$
V^{-1}f^+_2(y)=e^{i\pi/8}(\sqrt 2\pi h)^{-1/2}
\int_0^\infty e^{i\varphi(y,z)/h}dz,
\quad \varphi(y,z)=\frac 12(y^2+z^2-2\sqrt2yz+|\gamma|^2\ln z).
$$
Since $\gamma=O(h)$, the phase function
$\varphi(y,z)$ has two real critical points
$$
z_c^\pm(y)=(y\pm\sqrt{y^2-|\gamma|^2})/\sqrt2,
$$ 
which are the roots of $\frac{\partial\varphi}{\partial z}(y,z)= z-\sqrt2y+\frac{|\gamma|^2}{2z}$.
The phase of the asymptotic expansion will be given by
$\varphi(y,z_c^\pm(y))$.
Since we are on $\sigma_r^-$, we can assume $x>\sqrt E$, i.e. $y>0$ and independently from $h$. Then, 
$$
z_c^+(y)=\sqrt 2 y+O(h^2),\quad 
z_c^-(y)=\frac{|\gamma|^2}{2\sqrt2y}+O(h^2)
$$
and
$$
\varphi(y,z_c^\pm(y))=\mp\frac 12y^2+O(h^2).
$$
In view of Lemma \ref{bkw}, this means that the critical points
$z_c^+(y)$ and $z_c^-(y)$ will contribute on $\sigma_r^-$ and $\sigma_r^+$, respectively. 
Hence, we have only to compute the contribution from $z_c^+(y)$.
Since
$\frac{\partial^2\varphi}{\partial
z^2}(y,z_c^+(y))=1-\frac{|\gamma|^2}{2z_c^+(y)^2}=1+O(h^2)$,
$z_c^+(y)$ is a non-degenerate critical point, and 
the stationary phase theorem (e.g. Proposition 5.2 in \cite{dijo}) says that
$$
V^{-1}f^+_2(y)=e^{i\pi/8}\,2^{1/4}  \,
e^{i\varphi(y,z_c^+(y))/h}+O(h)
$$
microlocally near $\sigma_r^-$. Comparing $i\varphi(\phi(x),z_c^+(\phi(x)))=-\frac i2\phi(x)^2+O(h^2)$ with $z(x;r_2)=\int_{r_2}^x \sqrt{\nu^2-t^2(E-t^2)^2}/t \,dt$, we observe
for $x>r_2$
\begin{eqnarray*}
-\frac{d}{dx}z(x;r_2)
&=&
-\sqrt{\nu^2-x^2(E-x^2)^2}/x=-i(x^2-E)+O(h^2)\\
&=&
\frac{d}{dx}\left(-{\textstyle\frac i2}\phi(x)^2\right)+O(h^2),
\end{eqnarray*}
since $\nu^2=O(h^2)$ and $\phi(x)\phi'(x)=x^2-E$. Because of $z(r_2;r_2)=0$, one obtains 
$$
-z(x;r_2)=
-{\textstyle\frac i2}\phi(x)^2-{\textstyle\frac i2}\phi(r_2)^2 + O(h^2).
$$
If $\phi(r_2)^2=O(h^2)$, then $i\varphi(\phi(x),z_c^+(\phi(x)))=-z(x;r_2)+O(h^2)$. Consequently, 
$$
V^{-1}f^+_2(\phi(x))=e^{i\pi/8}\,2^{1/4}  \,
e^{-z(x;r_2)/h}\left(1+O(h)\right)
$$
and
$$
\tilde f^+(x)=e^{i\pi/8}\,
2^{1/4}  \, e^{-z(x;r_2)/h}
\begin{pmatrix}0\\ 1\end{pmatrix}
\left(1+O(h)\right)
$$
microlocally near $\sigma_r^-$, since $M(\phi(x),h)={\rm Id}+O(h)$.

\medskip
It remains to prove $\phi(r_2)^2=O(h^2)$. The turning point $r_2$ is a zero of the function $x\mapsto \nu^2/x^2 -(E-x^2)^2$. This is equivalent to ${\rm det}(a(r_2,0))=0$ with
$$
a(x,\xi)=\begin{pmatrix}
-\xi+x^2-E&\frac\nu x\\ -\frac\nu x& -\xi-x^2+E
\end{pmatrix}
$$
the symbol of the operator $A(x)-hD_x$. The first step of the normal form transformation of Theorem \ref{tono} reads on the symbol level as
$$
a(x,\xi)=\phi'(x)\tilde a(\phi(x),\xi/\phi'(x))
$$
with
$$
\tilde a(x,\xi)=\begin{pmatrix}
-\xi+x&\nu\psi(x)\\ -\nu\psi(x)& -\xi-x
\end{pmatrix}.
$$
Since $\phi'(r_2)\neq0$, one has ${\rm det}(a(r_2,0))=0$ if and only if ${\rm det}(\tilde a(\phi(r_2),0))=0$. The rest of the normal form transformation is
$$
\sqrt2\, q\circ\kappa_{\frac\pi4} =
M\,\sharp_h\, \tilde a \,\sharp_h\, M^{-1}
$$
with 
$$
q(y,\eta)=\begin{pmatrix}y&\frac\gamma{\sqrt2}\\ 
-\frac{\overline\gamma}{\sqrt2}& -\eta\end{pmatrix}
$$
the symbol of the normal form operator $Q$, and $\sharp_h$ the Moyal product of semiclassical Weyl calculus. The $h$-asymptotics of $\sharp_h$ (e.g. Proposition 7.7 in \cite{dijo}) together with the linear $\xi$-dependance of $\tilde a(x,\xi)$ yield
$$
M\,\sharp_h\, \tilde a\, \sharp_h\, M^{-1} =
M \tilde a  M^{-1} - i h M' M^{-1} = 
M \tilde a  M^{-1} + O(h^2).
$$
Hence, $\tilde a=\sqrt2\, M^{-1} (q\circ\kappa_{\frac\pi4}) M + O(h^2)$
and 
$$
{\rm det}(a(r_2,0))=0\;\Leftrightarrow\;
{\rm det}(q\circ\kappa_{\frac\pi4}(r_2,0))=O(h^2)
\;\Leftrightarrow\;
\phi(r_2)^2=O(h^2).
$$

%%%%%%%%%%%%%%%%%%%%%%%%%%%%%%%%%%%%%%%%%%%%%%%%%%%%%%%%%%%%%%%%%%%%%%%%%%%%%%%%%%%%%%%%%%%%%%%%%%%%%%%%%%%%%%%

%%%%%%%%%%%%%%%%%%%%%%%%%%%%%%%%%%%%%%%%%%%%%%%%%%%%%%%%%%%%%%%%%%%%%%%%%%%%%%

\end{document}